\documentclass[12pt]{article}
\oddsidemargin=\evensidemargin
\newfont{\blb}{msbm10 scaled\magstep1}

\newtheorem{theo}{Theorem}[section]

\newtheorem{prop}[theo]{Proposition}
\newtheorem{lemm}[theo]{Lemma}

\usepackage{amsmath}
\usepackage{amsfonts}
\usepackage{amssymb}

\date{}
\author{Anastasia Hadjievangelou, 
Patrizia Longobardi,  Mercede Maj, \\ Carmine Monetta,
E. A. O'Brien, and Gunnar Traustason 
}

\title{Left $3$-Engel elements in groups: A survey}

\begin{document}

\maketitle

\begin{abstract}{\bf Abstract.\ }We survey  left 3-Engel elements in groups.
\end{abstract}

\footnote{O'Brien is partially supported by the
Marsden Fund of New Zealand.}

\section{Introduction}

An element $a$ of a group $G$ is left Engel if for each $x\in G$ there exists a non-negative integer $n(x)$ such that
      $$[[[x,\underbrace {a],a],\ldots ,a]}_{n(x)}=1.$$
If $n(x)$ is bounded above, then $a$ is bounded left Engel. More precisely, $a$ is left $n$-Engel if $n(x)\leq n$ for all $x\in G$. \\ \\
Recall that the Hirsch-Plotkin radical, $\mbox{HP}(G)$, of $G$ is the product of its locally nilpotent normal subgroups. 
It is straightforward to see that every element of  $\mbox{HP}(G)$ is left Engel in $G$ and the converse is known
for some classes of groups, including solvable groups and 
finite groups (more generally, groups satisfying the maximal condition on
subgroups) \cite{Baer, Gru}. The converse is  not true in general, not even for bounded left Engel elements. In fact, whereas it is clear that a left $2$-Engel element is always in $\mbox{HP}(G)$, this remains an open question for left  $3$-Engel elements. There have been breakthroughs in recent years. In \cite{Jab} it is shown that every left $3$-Engel element of odd order is contained in $\mbox{HP}(G)$. In \cite{GA1}  this result is generalised to include elements of every order, by replacing left $3$-Engel with a stronger condition that we call strong left $3$-Engel. Left $3$-Engel elements of odd order are strong 
left $3$-Engel.   \\ \\
The results of \cite{Trac} imply that, in order to show that every left $3$-Engel element of finite order is in $\mbox{HP}(G)$, it suffices to consider elements of order $2$. By looking at a similar setting for Lie algebras, there are reasons to doubt that elements of order $2$ must be in $\mbox{HP}(G)$ however. Little is known about left $4$-Engel elements, although there are some interesting results in \cite{Ab2}.  \\ \\
Groups of prime power exponent satisfy some Engel type conditions and the solution to the restricted Burnside problem also uses the fact that the associated Lie ring satisfies certain Engel type identities \cite{zc,zd}. Burnside \cite{Burn} observed that every element of a group of exponent 3 is a left 2-Engel element, and so the fact that every left 2-Engel element lies in the Hirsch-Plotkin radical can be seen as the underlying reason why groups of exponent 3 are locally finite. For groups of 2-power exponent there is a close link with left Engel elements. If $a$ is an involution in a finitely generated group $G$ of exponent $2^{n}$, then 
            $$[[[x,\underbrace{a],a],\ldots ,a}_{n+1}]=[x,a]^{(-2)^{n}}=1.$$ 
Thus $a$ is a left $(n+1)$-Engel element of $G$. Hence,
if $G/G^{2^{n-1}}$ is finite and the left $(n+1)$-Engel elements of $G$ are
in $\mbox{HP}(G)$, then $G$ is finite. For sufficiently large $n$, the variety of groups of exponent $2^{n}$ is not locally finite \cite{Ivan,Lys}, so for sufficiently large $n$ there are left $n$-Engel elements not contained in $\mbox{HP}(G)$. Since groups of exponent 4 are locally finite \cite{San},  if all left 4-Engel elements of a group $G$ of exponent 8 are in $\mbox{HP}(G)$, then $G$ is locally finite. \\ \\
Is a left $3$-Engel element of prime power order in a group $G$ contained in $\mbox{HP}(G)$? We reduce this question to  elements of prime order using  the following result by Abdollahi \cite{Ab1}: for every prime $p$ and every left 3-Engel element $x$ of finite $p$-power, $x^p$ is in the Baer radical of $G$ and, in particular, $\langle x^p \rangle ^G$ is locally nilpotent. For the case $p=2$ this implies the local finiteness of groups of exponent 4, originally proved by Sanov \cite{San}. It is also proved in \cite{Ab1} that two left 3-Engel elements generate a nilpotent group of class at most 4, and this bound is achieved. \\ \\
We now swap the role of $a$ and $x$ in the definition of a left Engel element. Thus $a\in G$ is a {\it right Engel element} if for each $x\in G$ there exists a non-negative 
integer $n(x)$ such that 
    $$[a,_{n(x)} x]=1.$$
If $n(x)$ is bounded above by $n$, then $a$ is a right $n$-Engel element. By a result
of Heineken \cite{Hein1},  if $a$ is a right $n$-Engel element of $G$, then $a^{-1}$ is a left $(n+1)$-Engel
element. \\ \\
In \cite{New} Newell proved that if $a$ is a right $3$-Engel element of $G$, then $a\in \mbox{HP}(G)$;  in
fact he proved the stronger result that $\langle a\rangle^{G}$ is nilpotent of class at most $3$. A natural question arises whether $\langle a\rangle^{G}$ is nilpotent when $a$ is a left $3$-Engel element of $G$.  In \cite{GGM} it was shown that this is not the case by giving an example of a locally finite $2$-group with a left $3$-Engel element $a$ such that $\langle a\rangle^{G}$ is not nilpotent and in \cite{GA} this was generalised to an infinite family of examples. In \cite{GAM} this was extended to include every odd prime. Thus for each prime $p$ there is a locally finite $p$-group $G$ containing a left $3$-Engel element $a$ such that $\langle a\rangle^{G}$ is not nilpotent.\\ \\
The structure of this survey is as follows. In Section 2 we discuss (strong) sandwich elements in  Lie algebras and groups. We give a consistent power-conjugate presentation \cite[Section 9.4]{CGT-Handbook} of the largest 3-generator sandwich group, which is nilpotent of class 5 and at most 3 if all  elements have odd order. We then present one of our main results:  every finitely generated strong sandwich group is nilpotent. It follows that a strong left 3-Engel element of a group $G$ always lies in $\mbox{HP}(G)$; we mention applications for groups of exponents 5, 9 and 8. We do not know if all sandwich groups of rank 4 are nilpotent, but in Section 3 we give partial results. We show that there is a largest finite  sandwich group of rank $4$ generated by involutions. In Section 3.1 we obtain a new reduction theorem for sandwich groups of
rank $4$ with generators that have two commuting pairs. In Section 3.2 we report some results of computations in residually nilpotent sandwich groups of rank $5$ generated by involutions. In particular, we show that these are finite if there are at least $3$ commuting pairs among the generators. Finally, in Sections 4 and 5 we consider the global nilpotence problem for the normal closure of left 3-Engel elements in locally finite $p$-groups.

\section{Local nilpotence of (strong) sandwich groups and strong left 3-Engel elements}

The approach of   \cite{GA1, Jab}  is by working with  sandwich groups; these are group-theoretic analogues of  sandwich Lie algebras introduced by Kostrikin \cite{Ko1}. Kostrikin and Zel'manov \cite{Kos} proved that sandwich Lie algebras are locally nilpotent and this  is a key ingredient to both Kostrikin's solution to the restricted Burnside problem for groups of prime exponent \cite{Ko2}, and Zel'manov's general solution \cite{zc,zd}. 

\subsection{Sandwich Lie algebras}
 As for group commutators, we use the left normed convention for Lie products. \\ \\
{\bf Definition}. An element $a$ of a Lie algebra $L$ is  a {\it sandwich element} if  $axa=0$ and $axya=0$ for all $x,y \in L$.  A Lie algebra is a {\it sandwich Lie algebra} if it can be generated (as Lie algebra) by sandwich elements.\\ \\
%
The second condition is superfluous in odd characteristic:  $0=x(yaa)=xyaa-2xaya+xaay=2axya$. In characteristic $2$ it is needed as this example shows. \\  \\
{\bf Example}.  Consider the largest vector space $L= \langle a,b,c \rangle$ over ${\rm GF}(2)$ subject to $Id(c)$ being abelian, $bc=0$ and $bxb=axa=cxc=0$, for all $x \in L$. Then $L$ is a Lie algebra generated by  $a$, $b$, $ab$, $c(ab)^{n}$, and $c(ab)^{n}a$ for $\ n\geq 0$. To show that $L$ is non-nilpotent (or, equivalently, infinite dimensional), we give a concrete example where these generators are basis elements. Let $y=ab, u_n = c(ab)^n$ and $v_{n+1}=c(ab)^n a$. 
We have the following multiplication table for these basis elements of $L$. 
$$\begin{array}{lll}
    & a b=y, &  \\
    & ya = y b =0, & \\
    u_n a=v_{n+1}, & u_n  b = 0, & u_n  y = u_{n+1},\\ 
    v_n  a= 0, & v_n  b = u_n, & v_n  y = v_{n+1},\\
    u_n  v_m =0, & v_n  v_m = 0, & u_n  u_m = 0.\\
\end{array}
$$
%
It suffices to show that the Jacobi identity holds for these generators.  
Since the characteristic is 2, we deduce that:
$$\begin{array}{lllll}
    aby+bya+yab & = & 0, & & \\
    u_n ab +ab u_n+ bu_n a & =  & v_{n+1} b +y u_n  & =  & 2 u_{n+1}=0,\\
     v_n ab + ab v_n + b v_n a & =  & y v_n + u_n a & =  & 2u_{n+1}=0,\\
    u_n ay + ay u_n + y u_n a & =  & v_{n+1} y + u_{n+1} a & =  & 2v_{n+2}=0,\\
    v_n by+by v_n +y v_n b & = & u_n y+v_{n+1}b & = & 2u_{n+1}=0.
\end{array}
$$
%
\subsection{3-generator sandwich groups}

In \cite{Trau} the notion of a sandwich group was introduced. \\ \\
{\bf Definition}. A subset $X$  of a group $G$  is a {\it  sandwich set}  if  $\langle a,b^{g}\rangle$ is nilpotent of class at  most $2$ for all $a,b\in X$
and $g\in \langle X \rangle$. If $G$ is generated by a sandwich set, then $G$ is a {\it sandwich group}.  \\ \\
The connection with left $3$-Engel elements arises because the following are equivalent.   \vspace{2 mm} \\
(1) If $a$ is a left $3$-Engel element of a group $G$, then $a$ is in the locally nilpotent radical of $G$. \\
(2) Every finitely generated sandwich group is  nilpotent.\\ \\
The equivalence is a consequence of the following: if $a$ is a left $3$-Engel element in $G$, then $\langle a\rangle^{G}$ is a sandwich group; 
every element of a sandwich set $X$ is left $3$-Engel in $\langle X\rangle$. \\ \\
A significant feature of sandwich groups is that there is a largest sandwich group of any given rank \cite{Trau}. 
%
%

\begin{theo} Let $R = \langle x,y,z\rangle$ be the free sandwich group of rank $3$.
If the generators have odd order, then $R$ has class at most $3$, else $R$ has class $5$. 
\end{theo}
We obtained in \cite{Trau} the following consistent power-conjugate presentation for $R$. \\ \\
%
%
Let $e_1(z,z^x,y)=[z,x,y,y]$.\\ \\
\underline{Generators} \\ \\
$x_1=e_{1}(z,z^x,y),\quad  
x_2=e_{1}(x,x^y,z), \quad x_3=e_1(y,y^z,x)$, \\ \\
$x_4=[z,x,[z,y]], \quad x_5=[x,y,[x,z]], \quad x_6=[y,z,[y,x]], \\
 x_7=[z,x,y], \quad x_8=[z,y,x]$\\ \\
$x_9=[z,x], \quad x_{10}=[z,y], \quad x_{11}=[x,y]$ \\ \\
$x_{12}=x, \quad x_{13}=y, \quad x_{14}=z$. \\ \\
\underline{Relations} \\ \\
$x_1^2=x_2^2=x_3^2=1, \quad x_4^2=x_5^2=x_6^2=x_3x_2x_1$, \\ \\
$x_4^{x_{12}}=x_4x_3, \quad x_4^{x_{13}}=x_4x_1, \quad x_5^{x_{13}}=x_5x_1, \quad x_5^{x_{14}}=x_5x_2,\\
x_6^{x_{12}}=x_6x_3, \quad x_6^{x_{14}}=x_6x_2,$ \\ 
$x_7^{x_9}=x_7x_2x_3, \quad x_7^{x_{10}}=x_7x_1, \quad x_7^{x_{11}}=x_7x_1, \quad x_7^{x_{12}}=x_7x_5x_3x_2,\\ 
x_7^{x_{13}}=x_7x_1,\quad x_7^{x_{14}}=x_7x_4x_3x_2,$\\
$x_8^{x_9}=x_8x_3, \quad x_8^{x_{10}}=x_8x_2x_1, \quad x_8^{x_{11}}=x_8x_3, \quad x_8^{x_{12}}=x_8x_3,\\ 
x_8^{x_{13}}=x_8x_6x_3,\quad x_8^{x_{14}}=x_8x_4x_3,$\\
$x_9^{x_{10}}=x_9x_4, \quad x_9^{x_{11}}=x_9x_5, \quad x_9^{x_{13}}=x_9x_7, \quad x_{10}^{x_{11}}=x_{10}x_6, \quad x_{10}^{x_{12}}=x_{10}x_8,\\ 
x_{11}^{x_{14}}=x_{11}x_8x_7^{-1}x_6x_5x_4x_3, \quad x_{12}^{x_{13}}=x_{12}x_{11}, \quad x_{12}^{x_{14}}=x_{12}x_9^{-1}, \quad x_{13}^{x_{14}}=x_{13}x_{10}^{-1}.$\\ \\

\subsection{(Strong) sandwich groups and strong left 3-Engel elements}
 
The results of \cite{Jab} imply that a group $G$ generated by a finite sandwich set consisting of elements of odd order is nilpotent. Thus  every left $3$-Engel element of odd order in $G$ is in $\mbox{HP}(G)$. \\ \\ 
We wish to extend this result to include groups generated by elements of any order. We take a cue from the definition of sandwich algebras. In \cite{GA1} we introduced the notion of a strong sandwich set. \\ \\
{\bf Definition}. A subset $X$ of a group $G$  is a {\it strong sandwich set}  if:   \vspace{1 mm} \\
(1) $\langle a,b^{g}\rangle$ is nilpotent of class at most $2$ for all $a,b\in X$
and $g\in \langle X \rangle$;\\
(2) $\langle a,b^{f},c^{g}\rangle$ is nilpotent of class at most $3$ for all $a,b,c\in X$ and $f,g\in \langle X\rangle$.  \vspace{ 2 mm} \\
If $G$ is generated by a strong sandwich set, then it is a {\it strong sandwich group}. \\ \\
{\bf Remark}. A strong sandwich group is somewhat analogous to a sandwich Lie algebra:  if all the elements of $X$ have odd order, then condition (2) is superfluous. This follows from Theorem 2.1. \\ \\
There is also a connection to left $3$-Engel elements. We first need a definition.  \\ \\
{\bf Definition}. An element $a$ of a group $G$ is a {\it strong left $3$-Engel element} if: \vspace{2 mm} \\
(1) $\langle a,a^{g}\rangle$ is nilpotent of class at most $2$ for all  $g\in G$.\\
(2) $\langle a,a^{f},a^{g}\rangle$ is nilpotent of class at most $3$ for all $f,g\in G$.\\ \\
{\bf Remark}. Notice that $a$ is left $3$-Engel in $G$ if and only if it satisfies (1). Condition (2) is superfluous when $a$ has odd order. Therefore a left $3$-Engel element of odd order is a strong left $3$-Engel element. \\ \\
Strong left $3$-Engel elements and strong sandwich groups are related because the following are equivalent. \vspace{2 mm} \\
(1) If $a$ is a strong left $3$-Engel element of $G$, then $\langle a \rangle ^G$ is locally nilpotent. \\
(2) Every finitely generated strong sandwich group  is nilpotent.\\ \\
The following results obtained in \cite{GA1} about finitely generated strong sandwich groups generalise work of Jabara and Traustason \cite{Jab}.  
\begin{theo}
Every finitely generated strong sandwich group is nilpotent.
\end{theo}
\begin{theo}
If $a$ is a strong left $3$-Engel element of a group $G$, then $\langle a\rangle^{G}$ is locally nilpotent.
\end{theo}
A key ingredient in proving these results was the following \cite{GA1}.
\begin{prop} 
Let $X$ be a strong sandwich set in a group $G$ and let $a,b\in X$. Then $X\cup \{[a,b]\}$ is also a strong sandwich set in $G$. 
\end{prop}
We discuss briefly how Theorem 2.2 was proved using this proposition.
Let $X=\{x_{1}, x_{2},\ldots ,x_{r}\}$ be a strong sandwich set and let $\overline{X}$ consist of all commutators in $X$ (in any order and with any bracketing). Iterated use of Proposition 2.4 shows that $\overline{X}$ is a strong sandwich set. As every $3$-generator
strong sandwich group is nilpotent of class at most $3$, the Hall-Witt identity essentially reduces to the Jacobi identity. So for every $u,v,w\in \overline{X}$
   $$[u,[v,w]]=[u,v,w][u,w,v]^{-1}.$$ 
As a result, we applied the same basic approach as used  in Chanyshev's proof of the local nilpotence of sandwich Lie algebras \cite[Section 3.2]{Vaug}.
%
We mention applications to groups of exponents $5, 9$ and $8$. \\ \\
{\bf Theorem A} \cite{Trau}. {\it A group of exponent $5$ is locally finite if and only if it satisfies the law} 
     $$[z,[y,x,x,x],[y,x,x,x],[y,x,x,x]]=1.$$
{\bf Remark}. This result implies in particular that a group of exponent $5$ is locally finite if and only if all of its $3$-generator subgroups are finite; 
it was originally proved by Vaughan-Lee \cite{Vaug1}. \\ \\
{\bf Theorem B} \cite{Jab}. {\it Let $w$ be a word in $n$ variables $x_{1},\ldots ,x_{n}$ where the variety of groups satisfying the law $w^{3}=1$ is a locally finite variety of groups of exponent $9$. Then the same is true for the variety of groups on $n+1$ variables satisfying the law $(x_{n+1}^{3}w^{3})^{3}=1$.} \\ \\
{\bf Remark}. We can use Theorem B to construct an explicit sequence of words. Define the word $w_{n}=w_{n}(x_{1},\ldots , x_{n})$ in $n$ variables recursively by $w_{1}=x_{1}$
and $w_{n+1}=x_{n+1}^{3}w_{n}^{3}$. The variety of groups satisfying the law $x_{1}^{3}=1$ is locally finite by Burnside and by repeated application of Theorem B we see that, for each $n\geq 1$, the variety of groups satisfying the law $w_{n}^{3}$ is a locally finite variety of groups of exponent $9$. \\ \\
{\bf Theorem C} \cite{GA1}.  {\it Assume all groups of exponent $8$ satisfying a law $w=1$ in $n$ variables $x_{1},\ldots ,x_{n}$ are locally nilpotent.  Let $V$ be the  variety on $n+3$ variables satisfying both $[x_{n+1},w,w,w]=1$ and $\langle w,w^{x_{n+2}},w^{x_{n+3}}\rangle$  is nilpotent of class at most $3$. Then $V$ is locally nilpotent.} \\ \\
{\bf Remark}. Starting for example with $w=x_{1}^{4}$, Theorem C  gives us a sequence of locally nilpotent varieties.

\section{Left 3-Engel involutions}

Theorem 2.3 implies that if a left $3$-Engel element has odd order, then its normal closure is locally nilpotent. This observation  and \cite{Trac} imply that, to generalise this to left $3$-Engel elements of arbitrary finite order, it suffices to consider elements of order $2$. 
Despite Theorems 2.2 and 2.3, the question whether the normal closure of a left $3$-Engel involution is locally nilpotent remains open. The best general result is the following.
\begin{theo} Let $x$ be a left $3$-Engel involution in a group $G$. If $\langle x\rangle^{G}$ has no elements of order $8$, then $\langle x\rangle^{G}$ is locally nilpotent. 
\end{theo}
Even this apparently slight progress was surprisingly difficult to establish. The difficulty comes from the fact that we do not know initially that $\langle x\rangle^{G}$ is a $2$-group. Once we have a $2$-group, it is clear that $\langle x\rangle^{G}$ has exponent $4$ and thus it is locally finite by Sanov \cite{San}.  \\ \\
Consider sandwich groups generated by involutions. From \cite{Trau} we know that such a group is nilpotent if its rank is at most $3$. The largest sandwich group of rank $3$ generated by involutions has order $2^{13}$ and class $5$. 
\subsection{Rank 4 sandwich groups generated by involutions}
It remains unknown whether sandwich groups of rank $4$ generated by involutions are nilpotent. Nickel's nilpotent quotient algorithm \cite{NQ} is available in both {\sf GAP}  \cite{GAP}
and {\sc Magma} \cite{MAGMA}.  Using these implementations,
we proved that such a group  has a largest nilpotent quotient of order $2^{776}$ and class $13$. The calculation took 18 days of CPU time using {\sc Magma} 2.27-3 on a computer with a 2.6 GHz processor. \\ \\ 
We have some partial results on the nilpotence of such groups.
Before presenting them, we need the following definition. \\ \\
{\bf Definition}. Let $G$ be a sandwich group generated by a finite set
$X=\{ a_{1},\dots ,a_{r}\} $ of sandwich elements. The {\it commutativity graph}, $V(G)$, of $G$ is an (undirected) graph whose vertices are the generators $X$, and a pair of distinct vertices $a_{i}$ and 
$a_{j}$ are joined by an edge if and only if $a_{i}$ and $a_{j}$ commute. \\ \\
{\bf Remarks}. (1) The commutativity graph of the free $r$-generator sandwich
group has no edges. The largest $r$-generator sandwich group with commutativity graph the complete graph on $r$ vertices is the free abelian group of
rank $r$. \\ \\
(2) Let $H$ and $K$ be the largest $r$-generator sandwich groups with 
commutativity graphs $V(H)$ and $V(K)$ respectively. If $V(H)\subseteq V(K)$,
then $K$ is isomorphic to a quotient of $H$. \\ \\
Assume for the remainder of the section that $G$ is a sandwich group generated by $4$ involutions. If $V(G)$ is the complete graph, then $G$ is elementary abelian of order $16$. There is just one type of
commutativity graph with $5$ edges, namely

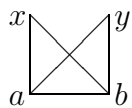
\mbox{}\\ \\
\begin{picture}(100,50)(-120,-50)
%
%
\put(15,-25){\line(1,0){30}}
\put(47,-29){$b$}
\put(7,-29){$a$}
\put(47,1){$y$}
\put(7,1){$x$}
\put(15,5){\line(1,-1){30}}
\put(15,-25){\line(1,1){30}}
\put(15,5){\line(0,-1){30}}
\put(45,5){\line(0,-1){30}}
\end{picture}
\normalsize
\mbox{}\\
and the largest sandwich group $\langle x,y,a,b\rangle$
with this commutativity graph is $\langle x,y\rangle\times \langle a,b\rangle=
D_{8}\times C_{2}^{2}$, a group of order $32$.  If the commutativity graph has $4$ edges, then there are two types of graphs to consider: either
the two removed edges are adjacent or not. \\
 \mbox{}\\ \\
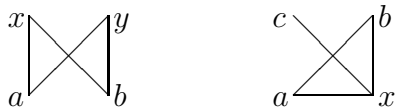
\begin{picture}(100,50)(-120,-50)
%
%
\put(27,-29){$b$}
\put(-13,-29){$a$}
\put(27,1){$y$}
\put(-13,1){$x$}
\put(-5,5){\line(1,-1){30}}
\put(-5,-25){\line(1,1){30}}
\put(-5,5){\line(0,-1){30}}
\put(25,5){\line(0,-1){30}}
%
\put(95,-25){\line(1,0){30}}
\put(127,-29){$x$}
\put(87,-29){$a$}
\put(127,1){$b$}
\put(87,1){$c$}
\put(95,5){\line(1,-1){30}}
\put(95,-25){\line(1,1){30}}
\put(125,5){\line(0,-1){30}}

\end{picture}
\normalsize
\mbox{}\\
The largest sandwich group with the first commutativity graph is $\langle x,y\rangle\times
\langle a,b\rangle=D_{8}\times D_{8}$, it has order $64$. For the second group, observe that $a$ and $b^{c}$ commute with each of $b$ and $a^{c}$. Thus 
     $$\langle a,b,c\rangle=\langle a,b^{c}\rangle \wr \langle c\rangle = D_{8}\wr C_{2},$$ 
 the standard wreath product of $D_{8}$ by $C_{2}$. Thus the largest sandwich group with the second commutativity graph is 
    $\langle a,b,c\rangle\times \langle x\rangle=(D_{8}\wr C_{2})\times C_{2},$
a group of order $256$. \\ \\
Next, we consider the case when the commutativity graph has 3 edges. This is 
much more difficult. There are three types
of commutativity graphs: \\
 \mbox{}\\ \\
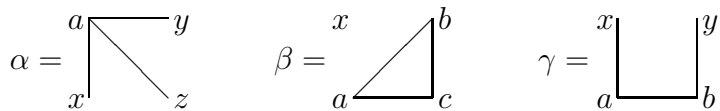
\begin{picture}(100,50)(-120,-50)
\put(-55,5){\line(1,0){30}}
\put(-23,-29){$z$}
\put(-63,-29){$x$}
\put(-23,1){$y$}
\put(-63,1){$a$}
\put(-85,-14){$\alpha=$}
\put(-55,5){\line(1,-1){30}}
\put(-55,5){\line(0,-1){30}}
%
\put(45,-25){\line(1,0){30}}
\put(77,-29){$c$}
\put(37,-29){$a$}
\put(15,-14){$\beta=$}
\put(77,1){$b$}
\put(37,1){$x$}
\put(45,-25){\line(1,1){30}}
\put(75,5){\line(0,-1){30}}
%
\put(145,-25){\line(1,0){30}}
\put(177,-29){$b$}
\put(137,-29){$a$}
\put(115,-14){$\gamma=$}
\put(177,1){$y$}
\put(137,1){$x$}
\put(145,5){\line(0,-1){30}}
\put(175,5){\line(0,-1){30}}
\end{picture}
\normalsize
\mbox{}\\
The largest sandwich group with the first commutativity graph is
    $$G_{\alpha}=\langle x,y,z\rangle\times \langle a\rangle =R\times C$$
where $R$ is the largest sandwich group generated by involutions $x,y,z$. From \cite{Trau} we know that $R$ is nilpotent of class at most $5$. The results of \cite{Trac} show that the others are finite. 
\subsubsection{Sandwich groups with commutativity graph $\beta$}
Let $G_{\beta}=\langle x,a,b,c\rangle$ be a sandwich group where $x,a,b,c$ are involutions and $G_{\beta}$ has commutativity graph
 \mbox{}\\ \\ \\
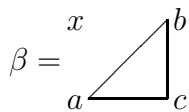
\begin{picture}(100,50)(-120,-50)
 \put(45,-25){\line(1,0){30}}
\put(77,-29){$c$}
\put(37,-29){$a$}
\put(15,-14){$\beta=$}
\put(77,1){$b$}
\put(37,1){$x$} 
\put(45,-25){\line(1,1){30}}
\put(75,5){\line(0,-1){30}}
\end{picture}
\normalsize
\mbox{}\\
In \cite{Trac} it was shown that $G_{\beta}$ has order $2^{28}$ and class $9$. It has the following consistent power-conjugate
presentation. \\ \\
Let 
     $$t(a)=[[x^{c},x^{a}],[x^{a},x^{b}]],\,t(b)=[[x^{a},x^{b}],[x^{b},x^{c}]],\,t(c)=[[x^{b},x^{c}],[x^{c},x^{a}]]$$
and
    $$y(a)=[x,x^{bc},x^{a}],\,y(b)=[x,x^{ca},x^{b}],\,y(c)=[x,x^{ab},x^{c}].$$
\mbox{}\\
\underline{Generators} \\ \\
$b_{1}=[t(a),b],\ 
b_{2}=t(a),\ b_{3}=t(b)$ \\
$b_{4}=y(b),\ b_{5}=y(b)^{a},\
b_{6}=y(a),\ b_{7}=y(a)^{b}, b_{8}=y(a)^{c},$ \\ \\
$b_{9}=[x,x^{ab}][x^{c},x^{abc}],\ b_{10}=[x,x^{ab}],\ 
b_{11}=[x,x^{bc}][x^{a},x^{abc}],$ \\
$b_{12}=[x,x^{bc}],\ b_{13}=[x,x^{ac}][x^{b},x^{abc}],\ 
b_{14}=[x,x^{ac}],\ b_{15}=[x,x^{abc}][x^{c},x^{ab}]$, \\
$b_{16}=[x,x^{abc}][x^{a},x^{bc}],\ 
b_{17}=[x,x^{abc}]$, \\ \\
$b_{18}=x,\ b_{19}=x^{a},\  b_{20}=x^{b},\ b_{21}=x^{c},\ 
b_{22}=x^{ab},\ b_{23}=x^{ca}$, \\
$b_{24}=x^{bc}, b_{25}=x^{abc}$, \\ \\
 $b_{26}=a,\ b_{27}=b,\ b_{28}=c$. \\ \\
\underline{Relations} \\ \\
$b_{1}^{2}=\ldots =b_{28}^{2}=1$. \\ \\
$b_{2}^{b_{27}}=b_{2}b_{1},\ b_{2}^{b_{28}}=b_{2}b_{1},\ 
b_{3}^{b_{26}}=b_{3}b_{1},\ b_{3}^{b_{28}}=b_{3}b_{1}$ \\
$b_{4}^{b_{19}}=b_{4}b_{2},\ b_{4}^{b_{21}}=b_{4}b_{3}b_{2},\ 
b_{4}^{b_{22}}=b_{4}b_{2}b_{1},
\ b_{4}^{b_{23}}=b_{4}b_{3},\ b_{4}^{b_{24}}=b_{4}b_{3}b_{2}b_{1}$, \\
$b_{4}^{b_{25}}=b_{4}b_{3},\ b_{4}^{b_{26}}=b_{5},\ b_{4}^{b_{28}}=b_{8}b_{6}b_{4}$ \\ 
$b_{5}^{b_{18}}=b_{5}b_{2},\ b_{5}^{b_{20}}=b_{5}b_{2}b_{1},\ b_{5}^{b_{21}}=
b_{5}b_{3}b_{1},\ b_{5}^{b_{23}}=b_{5}b_{3}b_{2}b_{1}$, \\
$b_{5}^{b_{24}}=b_{5}b_{3}b_{1},\ b_{5}^{b_{25}}=b_{5}b_{3}b_{2},\ 
b_{5}^{b_{26}}=b_{4},\ b_{5}^{b_{28}}=b_{8}b_{6}b_{5}$, \\
$b_{6}^{b_{20}}=b_{6}b_{3},\ b_{6}^{b_{21}}=b_{6}b_{3}b_{2},\ b_{6}^{b_{22}}=
b_{6}b_{3}b_{1},\ b_{6}^{b_{23}}=b_{6}b_{3}b_{2}b_{1},\ 
b_{6}^{b_{24}}=b_{6}b_{2}$, \\
$b_{6}^{b_{25}}=b_{6}b_{2},\  
b_{6}^{b_{27}}=b_{7},\ b_{6}^{b_{28}}=b_{8}$, \\
$b_{7}^{b_{18}}=b_{7}b_{3},\ b_{7}^{b_{19}}=b_{7}b_{3}b_{1},\ 
b_{7}^{b_{21}}=b_{7}b_{2}b_{1},\ 
b_{7}^{b_{23}}=b_{7}b_{2}b_{1}$ \\
$b_{7}^{b_{24}}=b_{7}b_{3}b_{2}b_{1},\ b_{7}^{b_{25}}=b_{7}b_{3}b_{2},\ 
b_{7}^{b_{27}}=b_{6},\ b_{7}^{b_{28}}=b_{8}b_{7}b_{6}$, \\
$b_{8}^{b_{18}}=b_{8}b_{3}b_{2},\ b_{8}^{b_{19}}=b_{8}b_{3}b_{2}b_{1},\ 
b_{8}^{b_{20}}=b_{8}b_{2}b_{1},\ b_{8}^{b_{22}}=b_{8}b_{2}b_{1}$, \\
$b_{8}^{b_{24}}=b_{8}b_{3}b_{1},\ b_{8}^{b_{25}}=b_{8}b_{3}, \
b_{8}^{b_{27}}=b_{8}b_{7}b_{6},\ 
b_{8}^{b_{28}}=b_{6}$, \\ \\
$b_{9}^{b_{12}}=b_{9}b_{1},\ 
b_{9}^{b_{14}}=b_{9}b_{1},\ b_{9}^{b_{18}}=b_{9}b_{6}b_{4},$ \\
$b_{9}^{b_{19}}=b_{9}b_{6}b_{5},\
b_{9}^{b_{20}}=b_{9}b_{7}b_{4},\
b_{9}^{b_{21}}=b_{9}b_{6}b_{4}$, \\
$b_{9}^{b_{22}}=b_{9}b_{7}b_{5},\ b_{9}^{b_{23}}=
b_{9}b_{6}b_{5},\ 
b_{9}^{b_{24}}=b_{9}b_{7}b_{4}$, \\
$b_{9}^{b_{25}}=b_{9}b_{7}b_{5}$, \\ 
$b_{10}^{b_{11}}=b_{10}b_{1},\ b_{10}^{b_{12}}=b_{10}b_{3},\ 
b_{10}^{b_{13}}=b_{10}b_{1},\ b_{10}^{b_{14}}=b_{10}b_{2}$, \\
$b_{10}^{b_{16}}=b_{10}b_{1},\ 
b_{10}^{b_{17}}=b_{10}b_{3}b_{2},\ b_{10}^{b_{21}}=b_{10}b_{6}b_{4}$, \\
$b_{10}^{b_{23}}=b_{10}b_{6}b_{5},\ 
b_{10}^{b_{24}}=b_{10}b_{7}b_{4},\ 
b_{10}^{b_{25}}=b_{10}b_{7}b_{5}$, \\
$b_{10}^{b_{28}}=b_{10}b_{9}$, \\
$b_{11}^{b_{14}}=b_{11}b_{1},\ b_{11}^{b_{18}}=b_{11}b_{6},\ 
b_{11}^{b_{19}}=b_{11}b_{6}$, \\
$b_{11}^{b_{20}}=b_{11}b_{7},\ b_{11}^{b_{21}}=b_{11}b_{8},\ 
b_{11}^{b_{22}}=b_{11}b_{7},\ b_{11}^{b_{23}}=b_{11}b_{8},\ 
b_{11}^{b_{24}}=b_{11}b_{8}b_{7}b_{6}$, \\
$b_{11}^{b_{25}}=b_{11}b_{8}b_{7}b_{6}$, \\
$b_{12}^{b_{13}}=b_{12}b_{1},\ b_{12}^{b_{14}}=b_{12}b_{3}b_{2},\ 
b_{12}^{b_{15}}=b_{12}b_{1}$, \\
$b_{12}^{b_{17}}=b_{12}b_{2},\ b_{12}^{b_{19}}=b_{12}b_{6},\ 
b_{12}^{b_{22}}=b_{12}b_{7},\ b_{12}^{b_{23}}=b_{12}b_{8},\ 
b_{12}^{b_{25}}=b_{12}b_{8}b_{7}b_{6}$, \\
$b_{12}^{b_{26}}=b_{12}b_{11}$, \\
$b_{13}^{b_{18}}=b_{13}b_{4},\ b_{13}^{b_{19}}=b_{13}b_{5},\ 
b_{13}^{b_{20}}=b_{13}b_{4}$ \\
$b_{13}^{b_{21}}=b_{13}b_{8}b_{6}b_{4},\ b_{13}^{b_{22}}=b_{13}b_{5},\ 
b_{13}^{b_{23}}=b_{13}b_{8}b_{6}b_{5},\ 
b_{13}^{b_{24}}=b_{13}b_{8}b_{6}b_{4}$ \\ 
$b_{13}^{b_{25}}=b_{13}b_{8}b_{6}b_{5}$, \\
$b_{14}^{b_{15}}=b_{14}b_{1},\ b_{14}^{b_{16}}=b_{14}b_{1},\ 
b_{14}^{b_{17}}=b_{14}b_{3},\ 
b_{14}^{b_{20}}=b_{14}b_{4}$, \\
$b_{14}^{b_{22}}=b_{14}b_{5},\ 
b_{14}^{b_{24}}=b_{14}b_{8}b_{6}b_{4},\ 
b_{14}^{b_{25}}=b_{14}b_{8}b_{6}b_{5},\ 
b_{14}^{b_{27}}=b_{14}b_{13}$, \\
$b_{15}^{b_{18}}=b_{15}b_{6}b_{4}b_{3}b_{2},\
b_{15}^{b_{19}}=b_{15}b_{6}b_{5}b_{3}b_{2}b_{1}$, \\
$b_{15}^{b_{20}}=b_{15}b_{7}b_{4}b_{3}b_{2}b_{1},\ 
b_{15}^{b_{21}}=b_{15}b_{6}b_{4}b_{3}b_{2},\ 
b_{15}^{b_{22}}=b_{15}b_{7}b_{5}b_{3}b_{2},\ 
b_{15}^{b_{23}}=b_{15}b_{6}b_{5}b_{3}b_{2}b_{1}$, \\ 
$b_{15}^{b_{24}}=b_{15}b_{7}b_{4}b_{3}b_{2}b_{1},\ 
b_{15}^{b_{25}}=b_{15}b_{7}b_{5}b_{3}b_{2}$, \\
$b_{16}^{b_{18}}=b_{16}b_{6}b_{2},\ 
b_{16}^{b_{19}}=b_{16}b_{6}b_{2},\ 
b_{16}^{b_{20}}=b_{16}b_{7}b_{2}b_{1}$, \\ 
$b_{16}^{b_{21}}=b_{16}b_{8}b_{2}b_{1},\ 
b_{16}^{b_{22}}=b_{16}b_{7}b_{2}b_{1},\ 
b_{16}^{b_{23}}=b_{16}b_{8}b_{2}b_{1},\ 
b_{16}^{b_{24}}=b_{16}b_{8}b_{7}b_{6}b_{2}$, \\ 
$b_{16}^{b_{25}}=b_{16}b_{8}b_{7}b_{6}b_{2}$, \\
$b_{17}^{b_{19}}=b_{17}b_{6}b_{2},\ 
b_{17}^{b_{20}}=b_{17}b_{4}b_{3},\ 
b_{17}^{b_{21}}=b_{17}b_{6}b_{4}b_{3}b_{2},\ 
b_{17}^{b_{22}}=b_{17}b_{7}b_{5}b_{3}b_{2}$, \\
$b_{17}^{b_{23}}=b_{17}b_{8}b_{6}b_{5}b_{3},\ 
b_{17}^{b_{24}}=b_{17}b_{8}b_{7}b_{6}b_{2},\ 
b_{17}^{b_{26}}=b_{17}b_{16},\ b_{17}^{b_{27}}=b_{17}b_{16}b_{15}$, \\
$b_{17}^{b_{28}}=b_{17}b_{15}$, \\ \\
$b_{18}^{b_{22}}=b_{18}b_{10},\ b_{18}^{b_{23}}=b_{18}b_{14},\ 
b_{18}^{b_{24}}=b_{18}b_{12},\ b_{18}^{b_{25}}=b_{18}b_{17},\ 
b_{18}^{b_{26}}=b_{19}$, \\
$b_{18}^{b_{27}}=b_{20},\ b_{18}^{b_{28}}=b_{21}$, \\
$b_{19}^{b_{20}}=b_{19}b_{10},\ b_{19}^{b_{21}}=b_{19}b_{14},\ 
b_{19}^{b_{24}}=b_{19}b_{17}b_{16},\ b_{19}^{b_{25}}=b_{19}b_{12}b_{11},\
b_{19}^{b_{26}}=b_{18}$, \\ 
$b_{19}^{b_{27}}=b_{22},\ b_{19}^{b_{28}}=b_{23}$, \\
$b_{20}^{b_{21}}=b_{20}b_{12},\ b_{20}^{b_{23}}=b_{20}b_{17}b_{16}b_{15},\ 
b_{20}^{b_{25}}=b_{20}b_{14}b_{13},\ b_{20}^{b_{26}}=b_{22},\ 
b_{20}^{b_{27}}=b_{18}$, \\ 
$b_{20}^{b_{28}}=b_{24}$, \\
$b_{21}^{b_{22}}=b_{21}b_{17}b_{15},\ b_{21}^{b_{25}}=b_{21}b_{10}b_{9},\
b_{21}^{b_{26}}=b_{23},\ b_{21}^{b_{27}}=b_{24},\ b_{21}^{b_{28}}=b_{18}$, \\
$b_{22}^{b_{23}}=b_{22}b_{12}b_{11},\ b_{22}^{b_{24}}=b_{22}b_{14}b_{13},\ 
b_{22}^{b_{26}}=b_{20},\ b_{22}^{b_{27}}=b_{19},\ 
b_{22}^{b_{28}}=b_{25}$, \\
$b_{23}^{b_{24}}=b_{23}b_{10}b_{9},\ b_{23}^{b_{26}}=b_{21},\ 
b_{23}^{b_{27}}=b_{25},\ b_{23}^{b_{28}}=b_{19}$, \\
$b_{24}^{b_{26}}=b_{25},\ b_{24}^{b_{27}}=b_{21},\ 
b_{24}^{b_{28}}=b_{20}$, \\
$b_{25}^{b_{26}}=b_{24},\ b_{25}^{b_{27}}=b_{23},\
b_{25}^{b_{28}}=b_{22}$. \\ 
\subsubsection{Sandwich groups with commutativity graph $\gamma$}
Let $G_{\gamma}=\langle a,b,x,y\rangle$ be a sandwich group where $a,b,x,y$ are involutions and $G_{\gamma}$
has  commutativity graph  \\ \\ \\
\begin{picture}(100,50)(-120,-50)
\put(45,-25){\line(1,0){30}}
\put(77,-29){$b$}
\put(37,-29){$a$}
\put(15,-14){$\gamma=$}
\put(77,1){$y$}
\put(37,1){$x$}
\put(45,5){\line(0,-1){30}}
\put(75,5){\line(0,-1){30}}
\end{picture}
\normalsize
\mbox{} \\ 
In \cite{Trac}, it was shown that $G_{\gamma}$ has order $2^{20}$ and class $9$. It  has the following consistent power-conjugate presentation.  \\ \\
\underline{Generators} \\ \\
$e_{1}=[x,b,[y,a],x,y,[x,b]],\ 
e_{2}=[x,b,[y,a],y,x,[y,a]]$, \\
$e_{3}=[x,b,[y,a],x,y,x],\ 
e_{4}=[x,b,[y,a],y,x,y]$, \\ 
$e_{5}=[x,b,[y,a],x,y],\ 
e_{6}=[x,b,[y,a],y,x]$, \\
$e_{7}=[x,b,[y,a],x],\ 
e_{8}=[x,b,[y,a],y]$, \\ 
$e_{9}=[x,b,y,x],\ e_{10}=
[x,[y,a],y]$ \\ \\
$e_{11}=[x,b,[y,a]],\ e_{12}=[x,b,y],\ 
e_{13}=[x,[y,a]]$, \\ 
$e_{14}=[x,y]$ \\ \\
$e_{15}=[x,b],\ e_{16}=[y,a],\ 
e_{17}=x,\ e_{18}=y$ \\ \\
$e_{19}=a,\ e_{20}=b$. \\ \\
\underline{Relations} \\ \\
$e_{1}^{2}=e_{2}^{2}=\ldots =e_{20}^{2}=1$. \\ \\
$e_{3}^{e_{20}}=e_{3}e_{1},\ e_{4}^{e_{19}}=e_{4}e_{2},\ 
e_{5}^{e_{15}}=e_{5}e_{1},\ e_{5}^{e_{17}}=e_{5}e_{3},$ \\
$e_{6}^{e_{16}}=e_{6}e_{2},\ e_{6}^{e_{18}}=e_{6}e_{4},\ 
e_{7}^{e_{12}}=e_{7}e_{1},\ e_{7}^{e_{14}}=e_{7}e_{3},\ 
e_{7}^{e_{18}}=e_{7}e_{5}$,\\
$e_{8}^{e_{13}}=e_{8}e_{2},\ e_{8}^{e_{14}}=e_{8}e_{4},\ e_{8}^{e_{17}}=e_{8}e_{6}$ \\
$e_{9}^{e_{11}}=e_{9}e_{1},\ e_{9}^{e_{13}}=e_{9}e_{3},\ 
e_{9}^{e_{16}}=e_{9}e_{5}e_{4}e_{2},\ e_{9}^{e_{19}}=
e_{9}e_{7}e_{6}e_{1}$, \\
$e_{10}^{e_{11}}=e_{10}e_{2},\ e_{10}^{e_{12}}=e_{10}e_{4},\
e_{10}^{e_{15}}=e_{10}e_{6}e_{3}e_{1},\ e_{10}^{e_{20}}=
e_{10}e_{8}e_{5}e_{2}$, \\
$e_{11}^{e_{14}}=e_{11}e_{6}e_{5}e_{4}e_{3},\ e_{11}^{e_{17}}=e_{11}e_{7},\ 
e_{11}^{e_{18}}=e_{11}e_{8}$, \\
$e_{12}^{e_{13}}=e_{12}e_{6}e_{5}e_{4}e_{3},\ 
e_{12}^{e_{16}}=e_{12}e_{8},\ e_{12}^{e_{17}}=e_{12}e_{9},\ 
e_{12}^{e_{19}}=e_{12}e_{11}e_{8}$, \\
$e_{13}^{e_{15}}=e_{13}e_{7},\ e_{13}^{e_{18}}=e_{13}e_{10},\ 
e_{13}^{e_{20}}=e_{13}e_{11}e_{7}$, \\
$e_{14}^{e_{15}}=e_{14}e_{9},\ e_{14}^{e_{16}}=e_{14}e_{10},\ 
e_{14}^{e_{19}}=e_{14}e_{13}e_{10},\ e_{14}^{e_{20}}=e_{14}e_{12}e_{9}$, \\
$e_{15}^{e_{16}}=e_{15}e_{11},\ e_{15}^{e_{18}}=e_{15}e_{12},\ 
e_{16}^{e_{17}}=e_{16}e_{13}$, \\
$e_{17}^{e_{18}}=e_{17}e_{14},\ e_{17}^{e_{20}}=e_{17}e_{15},\ 
e_{18}^{e_{19}}=e_{18}e_{16}$. 

\subsubsection{Sandwich groups whose commutativity graph have $2$ edges}
Let $G=\langle a,x,y,z\rangle$ be a sandwich group where $a,x,y,z$ are involutions and $G$ has 
 \mbox{}\\ \\ \\
\begin{picture}(100,50)(-120,-50)
%
%
\put(-23,-29){$y$}
\put(-63,-29){$a$}
\put(-23,1){$z$}
\put(-63,1){$x$}
\put(-85,-14){$\delta=$}
\put(-55,-25){\line(1,0){30}}
\put(-55,5){\line(0,-1){30}}
%
%
\end{picture}
\mbox{}\\
We know little about such groups.  In the remainder of this section we prove the following ``reduction" theorem. 
\normalsize
\begin{theo}  If $[[z,x],[z,a]]=1$ and $[[z,y],[z,a]]=1$, then $G$ is nilpotent. 
\end{theo}

\noindent 
{\bf Remark.} 
We used the {\sc Magma} implementation of the
nilpotent quotient algorithm to determine that
the largest nilpotent quotient $Q$ of 
$G$ has order $2^{71}$ and class $10$. 
If we impose the two relations of Theorem 3.2,
then $Q$ has order $2^{38}$ and class $9$. \\  

\noindent 
Let $N=\bigcup_{i=0}^{\infty}Z_{i}(G)$. As $G$ is finitely generated, it suffices to show that $G/N$ is nilpotent. Since $G$ is finitely generated, $G/N$ has trivial centre. Without loss of generality, we can thus assume that $Z(G)=1$.
Under this assumption, we will prove that $G=1$. Our lengthy proof of this claim concludes in Proposition 3.21 where we prove that $G$ is nilpotent. \\ 

The next two results do not rely on this assumption. The first is particularly useful. 
\begin{prop}
Let $u\in G$ where $u$ commutes with each of 
$x,y,[x,z],[y,z],$ and $[a,z]$. Then $u=1$.
\end{prop}
{\bf Proof}\ \ We first show that if $v\in G$ commutes with 
each of $x,y,z$ and $[z,a]$, then $v=1$. As a first step, we prove that  $[v,a]$ commutes with $x,y,z$ and $[z,a]$. Since each of $a$ and $v$ commutes with $x$, $y$ and $[z,a]$,  it is clear that $[v,a]$ commutes with each, so it suffices to show that
$[v,a]$ commutes with $z$. This follows from
               $$[v,a]^{z}=[v,a^{z}]=[v,a[a,z]]=[v,[z,a]a]=[v,a][v,[z,a]]^{a}=[v,a].$$
Thus each of  $v,[v,a]$ and $[v,_{2}a]$ commutes with each of  $x,y$ and $z$. 
However, $[v,_{2}a]$ commutes also with $a$.  Thus
$[v,_{2}a]\in Z(G)$ so $[v,_{2}a]=1$. It follows that $[v,a]=1$, so $v=1$. \\ \\ 
Next, notice that $[u,z]$ commutes with $x,y,[x,z],[y,z]$ and $[a,z]$. 
Note that $z$ commutes
with $[x,z], [y,z]$ and $[a,z]$, so $[u,z]$ commutes with each of
these elements. It thus suffices to show that $[u,z]$ commutes
with $x$ and $y$. This follows from
$$\begin{array}{l}
     [u,z]^{x}=[u,z[z,x]]=[u,[z,x]z]=[u,z][u,[z,x]]^{z}=[u,z], \\
\mbox{}
     [u,z]^{y}=[u,z[z,y]]=[u,[z,y]z]=[u,z][u,[z,y]]^{z}=[u,z].
\end{array}$$
Thus $[u,z,z]$ commutes with $x,y$ and $[z,a]$, but $[u,z,z]$ also commutes with $z$. From the first paragraph we deduce that 
 $[u,z,z]=1$. Hence $[u,z]$ commutes with each of $x,y,z$ and $[z,a]$ and thus it follows again
from the first paragraph that $[u,z]=1$. The same argument shows that $u=1$. $\Box$ \\ \\
Below we use the following power-conjugate presentation for the largest sandwich group of rank 3 generated by involutions; it is an immediate consequence of that given in Section 2.2. \\ \\
Let $e_1(z,z^x,y)=[z,x,y,y]$.\\ \\
\underline{Generators} \\ \\
$x_1=e_{1}(z,z^x,y),\quad  
x_2=e_{1}(x,x^y,z), \quad x_3=e_1(y,y^z,x)$, \\ \\
$x_4=[z,x,[z,y]], \quad x_5=[x,y,[x,z]], \quad x_6=[y,z,[y,x]], \\
x_7=[z,x,y], \quad x_8=[z,y,x]$\\ \\
$x_9=[z,x], \quad x_{10}=[z,y], \quad x_{11}=[x,y]$ \\ \\
$x_{12}=x, \quad x_{13}=y, \quad x_{14}=z$. \\ \\
\underline{Relations} \\ \\
$x_1^2=\dotsm =x_6^2=1, \quad x_3=x_2x_1, \quad 
x_7^2=x_1, \quad x_8^2=x_3, \quad x_9^2=\dotsm = x_{14}^2=1,$\\
$x_4^{x_{12}}=x_4x_2x_1, \quad x_4^{x_{13}}=x_4x_1, \quad x_5^{x_{13}}=x_5x_1, \quad x_5^{x_{14}}=x_5x_2,\\
x_6^{x_{12}}=x_6x_2x_1, \quad x_6^{x_{14}}=x_6x_2,$ \\ 
$x_7^{x_9}=x_7x_1, \quad x_7^{x_{10}}=x_7x_1, \quad x_7^{x_{11}}=x_7x_1, \quad x_7^{x_{12}}=x_7x_5x_1,\\ 
x_7^{x_{13}}=x_7x_1,\quad x_7^{x_{14}}=x_7x_4x_1,$\\
$x_8^{x_9}=x_8x_2x_1, \quad x_8^{x_{10}}=x_8x_2x_1, \quad x_8^{x_{11}}=x_8x_2x_1, \quad x_8^{x_{12}}=x_8x_2x_1,\\ 
x_8^{x_{13}}=x_8x_6x_2x_1,\quad x_8^{x_{14}}=x_8x_4x_2x_1,$\\
$x_9^{x_{10}}=x_9x_4, \quad x_9^{x_{11}}=x_9x_5, \quad x_9^{x_{13}}=x_9x_7, \quad x_{10}^{x_{11}}=x_{10}x_6, \quad x_{10}^{x_{12}}=x_{10}x_8,\\ 
x_{11}^{x_{14}}=x_{11}x_8x_7x_6x_5x_4x_2, \quad x_{12}^{x_{13}}=x_{12}x_{11}, \quad x_{12}^{x_{14}}=x_{12}x_9, \quad x_{13}^{x_{14}}=x_{13}x_{10}.$ \\ \\ 

\begin{lemm} $\gamma_{3}(\langle [x,y],[z,a]\rangle)=1$.  
\end{lemm}
{\bf Proof}\ \ Notice that $1=[x,y^{2}]=[x,y]^{2}$ and $1=[z,a^{2}]=[z,a]^{2}$. So
            $$[[x,y],[a,z],[a,z]]=[[x,y],[z,a],[x,y]]=([a,z][x,y])^{4}$$
and thus it suffices to show that $[[x,y],[z,a],[z,a]]\in Z_{2}(G)=1$. Calculating in the sandwich group $\langle
x,y,a^{z}\rangle$, using the presentation above, we see that  
     $$[[x,y],[a,z],[a,z]]=[x,y,a^{z},a^{z}]=[x_{11},x_{14},x_{14}]=[x_{8}x_{7}x_{6}x_{5}x_{4}x_{2},x_{14}]=x_{2}.$$
Observe that $x_{2}$ commutes with $x,y$ and $a^{z}$, and also with $a$. We finish the proof by showing that $[[x,y],[a,z],[a,z],z]$ commutes with 
each of $a,x,y$ and $z$. We will calculate in the 
sandwich group $\langle \bar{a},\bar{b},\bar{c},\bar{x}\rangle$,
where $\bar{a}=a,\ \bar{b}=x,\ \bar{c}=x^{y}$ and $\bar{x}=z$. Since it has commutativity graph $\beta$,
we can use the relations satisfied by such a group. We have already seen that 
$[[x,y],[a,z],[a,z]]$ commutes with $[z,a]=[\bar{x},\bar{a}]$ and $x=\bar{b}$.
We will use this in the following calculation:
\begin{eqnarray*}
 \mbox{}[[x,y],[z,a],[z,a]] & = & [\bar{b}\bar{c},[\bar{x},\bar{a}],
               [\bar{x},\bar{a}]] \\
           & = & [[\bar{x},\bar{a}]^{\bar{b}\bar{c}},
                  [\bar{x},\bar{a}]] \mbox{\ \ \ (commutes with $\bar{b}=x$)} \\
     & = & [[\bar{x},\bar{a}]^{\bar{c}},[\bar{x},\bar{a}]^{\bar{b}}] \\
   & = & [\bar{x}^{\bar{c}}\bar{x}^{\bar{a}\bar{c}},\bar{x}^{\bar{b}}
           \bar{x}^{\bar{a}\bar{b}}] \\
  & = & [\bar{x}^{\bar{c}},\bar{x}^{\bar{b}}x^{\bar{a}\bar{b}}]^{\bar{x}^
  {\bar{a}\bar{c}}}
[\bar{x}^{\bar{a}\bar{c}},\bar{x}^{\bar{b}}\bar{x}^{\bar{a}
\bar{b}}] \\
   & = & [\bar{x}^{\bar{c}},\bar{x}^{\bar{a}\bar{b}}]^
   {\bar{x}^{\bar{a}\bar{c}}}[\bar{x}^{\bar{c}},\bar{x}^{\bar{b}}]
    ^{\bar{x}^{\bar{a}\bar{b}}\bar{x}^{\bar{a}\bar{c}}}[\bar{x}^{\bar{a}\bar{c}},
\bar{x}^{\bar{a}\bar{b}}][\bar{x}^{\bar{a}\bar{c}},\bar{x}^{\bar{b}}]^{\bar{x}^
  {\bar{a}\bar{b}}} \\
  & = & (b_{17}b_{15})^{b_{23}}b_{12}^{b_{22}b_{23}}b_{12}^{b_{26}}
     (b_{17}b_{16}b_{15})^{b_{22}} \\
   & = & b_{16}b_{11}b_{2}b_{1}.
\end{eqnarray*}
This element commutes with $\bar{b}=x$ and thus 
      $$b_{16}b_{11}b_{2}b_{1}=(b_{16}b_{11}b_{2}b_{1})^{b_{27}}=
     b_{16}b_{11}b_{2}b_{1}b_{1},$$
so $b_{1}=1$. Hence 
    $$[[x,y],[z,a],[z,a],z]=[b_{16}b_{11}b_{2},b_{18}]=b_{2}$$
commutes with each of $x,z$ and $a$. Since $[x,y]=[y,x]$, 
we see by symmetry that $b_{2}$ 
also commutes with $y$. $\Box$ \\ \\
We deduce that $b_{2}=1$. We analyse some consequences. Consider the elements:
$$\begin{array}{l}
   b_{11}=[\bar{x}^{\bar{b}\bar{c}},\bar{x}][\bar{x}^{\bar{a}},\bar{x}^{\bar{a}\bar{b}\bar{c}}]=[z,z^{[x,y]}][z^{a},z^{a[x,y]}];  \\
\mbox{}b_{6}=[b_{11},b_{18}]=[b_{11},\bar{x}]=[b_{11},z]; \\
\mbox{}b_{8}b_{7}=[b_{6},b_{27}b_{28}]=[b_{6},[x,y]];
\end{array}$$
all are symmetric in $x,y$. \\ \\
Since $b_{1}=b_{2}=1$, we deduce that  $b_{8}b_{7}$ commutes with  each of $x,a$ and
$z$, and by symmetry also with $y$; thus $b_{8}b_{7}=1$. Recall that $1=[x,y,a^{z},a^{z}]=x_{2}$ in $\langle x,y,a^{z}\rangle$. In summary:
\begin{lemm} \mbox{}\\ 
{\rm (1)} $x_{2}=1$ in $\langle x,y,a^{z}\rangle $. \\
{\rm (2)} $b_{1}=b_{2}=1$, $b_{8}=b_{7}$ and $b_{16}=b_{11}$ in $\langle a,x,x^{y},z\rangle$ and $\langle a,y,y^{x},z\rangle$.
\end{lemm} 
From now on, we use the conditions given in the statement of Theorem 3.2. Since $[z^{y},z^{a}]=[z,y,[z,a]]=1$, we deduce that $\langle z, z^{a}, z^{y},x\rangle$ is a $\beta$-group. We will later work with this subgroup. First we prove a crucial
lemma. 
\begin{lemm} $[z,a,x]=[z,x,a]$ and $[z,a,y]=[z,y,a]$.
\end{lemm}
{\bf Proof} Since $[z,x,[z,a]]=1$, it follows that $[z,x,a]=[x^{z},a]$ commutes with $z$. Hence
                         $$[z,x,a]=[x^{z},a]=[x^{z},a]^{z}=[x,a^{z}]=[a^{z},x]=[z,a,x].$$
The latter claim is proved similarly. $\Box$
\begin{lemm} $[x,y,z,z]=1$.
\end{lemm}
{\bf Proof}\ \ We calculate in the $\beta$-group $\langle \bar{a}, \bar{b}, \bar{c}, \bar{x}\rangle=\langle a,x,x^{y},z\rangle$. We know that 
                        $$b_{1}=b_{2}=1,\ b_{8}=b_{7},\ b_{16}=b_{11}.$$
Furthermore 
      $$1=[z^{a},z^{x}]=[\bar{x}^{\bar{a}},\bar{x}^{\bar{b}}]=b_{10}.$$
As a consequence,
\begin{eqnarray*}
      1=[b_{10},b_{12}] & = & b_{3} \\
     1=[b_{10},b_{28}] & = & b_{9} \\
     1=[b_{9},b_{19}] & = & b_{6}b_{5} \\
     1=[b_{9},b_{22}] & = & b_{7}b_{5} \\
     1=[b_{9},b_{20}] & = & b_{7}b_{4}.
\end{eqnarray*}
The last three identities imply that $b_{4}=b_{5}=b_{6}=b_{7}=b_{8}$. We show that this common element is trivial. We consider  
                         $$[z,z^{[x,y]},z^{a}]=[\bar{x},\bar{x}^{\bar{b}\bar{c}},\bar{x}^{\bar{a}}]=[b_{12},b_{19}]=b_{6}.$$
From the presentation, $b_{6}^{x}=b_{6}^{\bar{b}}=b_{6}^{b_{27}}=b_{7}=b_{6}$, 
so $b_{6}$ commutes with $x$. Also $b_{6}$ commutes with $[z,x]=b_{18}b_{20}$ and 
$[z,a]=b_{18}b_{19}$. As $b_{6}$ is symmetric in $x$ and $y$, it commutes with $y$ and $[z,y]$.  Proposition 3.3 implies that $b_{6}=1$. Hence 
                         $$b_{4}=b_{5}=b_{6}=b_{7}=b_{8}=1.$$
Next, consider $[z,z^{[x,y]}][z^{a},z^{a[x,y]}]=b_{12}b_{12}^{b_{26}}=b_{11}$. Observe that $b_{11}$ commutes with each of $x=b_{27}$, $[z,x]=b_{18}b_{20}$ and $[z,a]=b_{18}b_{19}$. As $b_{11}$ is symmetric
in $x$ and $y$, it commutes with $y$ and $[z,y]$. Proposition 3.3 implies that $b_{11}=1$. Hence
                               $$b_{16}=b_{11}=1.$$
Finally, we consider $[x,y,z,z]=[z,z^{[x,y]}]=b_{12}$. Observe that $b_{12}$ commutes with each of $x=b_{27}$, $[z,x]=b_{18}b_{20}$ and $[z,a]=b_{18}b_{19}$. 
As $b_{12}$ is symmetric in $x$ and $y$, it commutes with $y$ and $[z,y]$. Proposition 3.3 implies that $b_{12}=[x,y,z,z]=1$. $\Box$ \\ \\
{\bf Remark}. We know from the presentation for $\langle x,y,z\rangle$ that $\gamma_{5}(\langle x,y,z\rangle)$ is generated by $[x,y,z,z], [z,x,y,y]$ and $[z,y,x,x]$, and their product is trivial. Also $[x,y,z,z]=1$. Hence 
$\gamma_{5}(\langle x,y,z\rangle)$ is cyclic, and generated by $[z,x,y,y]=[z,y,x,x]$. 
\begin{lemm} $\langle x,y,a^{z}\rangle$ is nilpotent of class at most $4$ and $[[a^{z},x],[a^{z},y]]=1$. 
\end{lemm}
{\bf Proof}\ \ As $\langle x,y,a^{z}\rangle$ is a sandwich group, we know from the presentation of the largest such group that $\gamma_{5}(\langle x,y,a^{z}\rangle)$ is generated by
$[x,y,a^{z},a^{z}],[a^{z},x,y,y]$
and $[a^{z},y,x,x]$. The product of these elements is trivial and 
Lemma 3.5 implies that $[x,y,a^{z},a^{z}]=1$. Using Lemma 3.6 and the presentation for the $\gamma$-group $\langle \bar{x},\bar{a},\bar{b},\bar{y}\rangle=\langle x^{z},x,a,y\rangle$, we deduce that 
           $$[a^{z},x,y,y]=[x^{z},a,y,y]=[\bar{x},\bar{b},\bar{y},\bar{y}]=[e_{12},e_{18}]=1.$$
Hence $\langle x,y,a^{z}\rangle$ is nilpotent of class at most $4$. For the last claim,
 notice that $[a^{z},x,[a^{z},y]]$ commutes with $z$ by the assumptions of Theorem 3.2 and, as $\langle x,y,a^{z}\rangle$ is nilpotent of class at most $4$, 
it commutes with $x$ and $y$, and it obviously commutes with $a$. Hence $[a^{z},x,[a^{z},y]]\in Z(G)=1$. $\Box$

\begin{lemm}
In the $\beta$-group $\langle \bar{a},\bar{b},\bar{c}, x\rangle =\langle z,z^{a},z^{y},x\rangle$, the following relations hold:
                    $$b_{1}=b_{2}=b_{3}=b_{4}=b_{5}=b_{6}=b_{7}=b_{9}=b_{10}=1.$$
\end{lemm}
{\bf Proof}\ \ Observe that $1=[z^{x},z^{a}]=[\bar{a}^{x},\bar{b}]=[x,\bar{a},\bar{b}]=xx^{\bar{a}}x^{\bar{b}}x^{\bar{a}\bar{b}}$ from which it follows that $x^{\bar{a}\bar{b}}=xx^{\bar{a}}x^{\bar{b}}$. In particular, 
                $$1=[x,x^{\bar{a}\bar{b}}]=b_{10}.$$
As a consequence,
\begin{eqnarray*}
     1=[b_{10},b_{11}] & = & b_{1} \\
     1=[b_{10},b_{14}] & = & b_{2} \\
     1=[b_{10},b_{12}] & = & b_{3} \\
     1=[b_{10},b_{28}] & = & b_{9} \\
     1=[b_{9},b_{18}] & = & b_{6}b_{4} \\ 
     1=[b_{9},b_{19}] & = & b_{6}b_{5} \\ 
     1=[b_{9},b_{20}] & = & b_{7}b_{4}.
\end{eqnarray*}
From the last three identities, $b_{4}=b_{5}=b_{6}=b_{7}$. We show that this common element is trivial. We consider $t=[z,a,x,[z,y,x,x]]$. Calculating first in $\langle z,z^{a},z^{y},x\rangle$, 
                  $$t=[\bar{a}\bar{b},x,[\bar{a}\bar{c},x,x]]=[b_{26}b_{27},b_{18},b_{14}]=[b_{18}^{b_{26}b_{27}}b_{18},b_{14}]=[b_{22}b_{18},b_{14}]=b_{5}.$$
From the presentation, $[b_{5},[z,a]]=[b_{5},b_{26}b_{27}]=b_{4}b_{5}=1$. \\ \\
We next consider $t$ in the $\gamma$-group $\langle \bar{x},\bar{a},\bar{b},\bar{y}\rangle=\langle x^{z},x,a,y\rangle$. We first find some relations in this group. Using the fact that $\langle x,y,a^{z}\rangle$ is nilpotent of class at most $4$ and Lemma 3.6, we deduce that 
        $$1=[a^{z},x,[x,y],y]=[x^{z},a,[x,y],y]=[\bar{x},\bar{b},[\bar{a},\bar{y}],\bar{y}]=[e_{11},e_{18}]=e_{8}.$$
As a consequence, 
\begin{eqnarray*}
                              1=[e_{8},e_{13}] & = & e_{2} \\
                              1=[e_{8},e_{14}] & = & e_{4} \\
                              1=[e_{8},e_{17}] & = & e_{6}.
\end{eqnarray*}
We use the facts that $[z,a,x]=[z,x,a]$ from Lemma 3.6 and $[z,y,x,x]=[z,x,y,y]$. Hence
\begin{eqnarray*}
            t & = & [[z,x,a],[z,x,y,y]] \\
              & = & [[x^{z},a], [x^{z}x,y,y]] \\
             & = & [[\bar{x},\bar{b}],[\bar{a}\bar{x},\bar{y},\bar{y}]] \\
             & = & [\bar{x},\bar{b}, [[\bar{a},\bar{y}][\bar{a},\bar{y},\bar{x}][\bar{x},\bar{y}],\bar{y}]] \\
             & = & [\bar{x},\bar{b},[\bar{a},\bar{y},\bar{x},\bar{y}]^{[\bar{x},\bar{y}]}] \\
             & = & [e_{15}, e_{10}^{e_{14}}] \\
             & = & [e_{15},e_{10}] \\
            & = & e_{3}e_{1}.
\end{eqnarray*}
From the presentation, we deduce that $t$ commutes with $y=\bar{y}=e_{18}$ and $x=\bar{a}=e_{19}$. Both $[z,x,a]$ and $[z,x,y,y]$ commute with $z$. Hence $t$ commutes with $z$ and thus with $[z,x]$ and $[z,y]$. Thus
$t$ commutes with each of  $x,y,[z,x],[z,y]$ and $[z,a]$. Proposition 3.3 implies that $t=1$. Hence $b_{5}=1$. $\Box$
\begin{lemm} $\langle x,y,z\rangle$ is nilpotent of class at most $4$ and $[[z,x],[z,y]]=1$. 
\end{lemm}
{\bf Proof}\ \ In the $\beta$-group $\langle \bar{a},\bar{b}, \bar{c}, x\rangle=\langle z,z^{a},z^{y},x\rangle$ 
                      $$[z,y,x,x]=[x^{zz^{y}},x]=[\bar{x}^{\bar{a}\bar{c}},x]=b_{14}.$$
The expressions for both $b_{14}$ and  $b_{13}=[z,y,x,x,[z,a]]=[b_{14},b_{26}b_{27}]$ are symmetric in $x$ and $y$. We read  from the presentation that:
\begin{eqnarray*}
      [b_{13},[z,a]] & = & [b_{13},b_{26}b_{27}]=1; \\
   \mbox{} [b_{13},x] & = & [b_{13},b_{18}]=b_{4}=1; \\
 \mbox{} [b_{13},[z,x]] & = & [b_{13},[\bar{a}, x]]=[b_{13},b_{18}b_{19}]=b_{5}b_{4}=1.
\end{eqnarray*}
As $b_{13}$ is symmetric in $x$ and $y$, it commutes with $y$ and $[z,y]$,
so by Proposition 3.3 it is trivial. 
Hence $[z,y,x,x]=b_{14}$ commutes with $[z,a]$. It clearly commutes with 
each of $x,y,[z,x]$ and $[z,y]$, so 
by Proposition 3.3 it is trivial. 
Therefore $\langle x,y,z\rangle$ is nilpotent of class at most $4$. 
%
Hence $[z,x,[z,y]]$ commutes with each of $x,y,[z,x]$ and $[z,y]$. As $[z,a]$ commutes with both $[z,x]$ and $[z,y]$ it follows that $[z,x,[z,y]]$ commutes with
$[z,a]$. Proposition 3.3 implies that $[z,x,[z,y]]=1$. $\Box$
\begin{lemm} In the $\gamma$-group $\langle \bar{x},\bar{a},\bar{b},\bar{y}\rangle=\langle x^{z},x,a,y\rangle$  the following relations hold:
  $$e_{1}=e_{2}=e_{3}=e_{4}=e_{5}=e_{6}=e_{7}=e_{8}=e_{10}=1,\ e_{11}=e_{9}.$$
\end{lemm}
{\bf Proof}\ \ We have seen in the proof of Lemma 3.9 that $e_{2}=e_{4}=e_{6}=e_{8}=1$. 
As $\langle x,y,z\rangle$ is nilpotent of class at most $4$, 
 \begin{eqnarray*}
        1 & = & [z,x,y,y] \\
           & = & [x^{z}x,y,y] \\
           & = & [\bar{x}\bar{a},\bar{y},\bar{y}] \\
           & = & [[\bar{x},\bar{y}][\bar{x},\bar{y},\bar{a}][\bar{a},\bar{y}],y] \\
          & = & [\bar{x},\bar{y},\bar{a},\bar{y}]^{[\bar{a},\bar{y}]} \\
          & = & [e_{14},e_{19},e_{18}]^{e_{16}} \\
         & = & [e_{13}e_{10},e_{18}]^{e_{16}} \\
         & = & e_{10}^{e_{16}} \\
        & = & e_{10}.
\end{eqnarray*}
Also 
$$\begin{array}{l}
       1=[e_{10},e_{15},e_{20}]  = [e_{3}e_{1},e_{20}]=e_{1} \\
      1=[e_{10},e_{15}]  =  e_{3} \\
      1=[e_{10},e_{20}] = e_{5}.
\end{array}$$
We consider $[z,x,a,[z,x,[x,y]]]$. Both $[z,x,a]$ and $[z,x,[x,y]]$ commute with each of $z$ and $x$ so
$[z,x,a,[z,x,[x,y]]]$ commutes with each. We compute in the $\gamma$-group $\langle x^{z},x,a,y\rangle$. Note that 
    $$[z,x,a,[z,x,[x,y]]]=[x^{z},a,[x^{z},[x,y]]]=[\bar{x},\bar{b},[\bar{x},[\bar{a},\bar{y}]]]=[e_{15},e_{13}]=e_{7}.$$
From the presentation, $e_{7}$ commutes with $e_{20}=\bar{b}=a$ and $e_{18}=\bar{y}=y$. It follows that
$e_{7}\in Z(G)=1$. \\ \\
Finally, we consider $[z,x,a,[z,x,y]]$. By Lemma 3.10 we know that $[z,x,y]$ commutes with $z$ and as $[z,x,a]$ also commutes with $z$ we see that $[z,x,a,[z,x,y]]$ commutes with $z$. From the $\gamma$-group presentation, we deduce that 
 \begin{eqnarray*}
  [z,x,a,[z,x,y]] & = & [x^{z},a,[x^{z}x,y]] \\
                      & = & [\bar{x},\bar{b}, [\bar{x}\bar{a},\bar{y}]] \\
                       & = & [e_{15},e_{14}^{e_{19}}e_{16}] \\
                       & = & [e_{15},e_{14}e_{13}e_{16}] \\
                      & = & e_{9}e_{11}.
\end{eqnarray*}
%
We also see that $e_{9}e_{11}$ commutes with 
$e_{20}=\bar{b}=a$, $e_{18}=\bar{y}=y$ and
$e_{17}=x$. Hence $e_{9}e_{11}\in Z(G)=1$. $\Box$
\begin{lemm} In the $\beta$-group $\langle \bar{a},\bar{b},\bar{c},\bar{x}\rangle=
\langle a,x,x^{y},z\rangle$ the following relations hold:
   $$b_{1}=b_{2}=b_{3}=b_{4}=b_{5}=b_{6}=b_{7}=b_{8}=b_{9}=b_{10}=b_{11}=b_{12}=b_{16}=1,\ b_{15}=b_{13}.$$
\end{lemm}
{\bf Proof}\ \  All equalities were  established in the proof of Lemma 3.7 apart from $b_{15}=b_{13}$. 
In the proof of Lemma 3.11 we show that $[z,x,a,[z,x,[x,y]]]=e_{7}=1$. We now view it
as an element of $\langle a,x,x^{y},z\rangle$. Our relations already imply that $\gamma_{3}(\langle \bar{x}\rangle^{\langle \bar{a},\bar{b},\bar{c}\rangle})=1$. This simplifies our computations: 
 \begin{eqnarray*}
   1 & = & [z,x,a,[z,x,[x,y]]] \\
     & = & [\bar{x}\bar{x}^{\bar{a}}\bar{x}^{\bar{b}}\bar{x}^{\bar{a}\bar{b}},
                \bar{x}\bar{x}^{\bar{b}}\bar{x}^{\bar{b}\bar{c}}\bar{x}^{\bar{c}}] \\
    & = & [\bar{x}^{\bar{a}},\bar{x}^{\bar{b}\bar{c}}][\bar{x}^{\bar{a}},\bar{x}^{\bar{c}}]
              [\bar{x}^{\bar{a}\bar{b}},\bar{x}^{\bar{b}\bar{c}}][\bar{x}^{\bar{a}\bar{b}},\bar{x}^{\bar{c}}] \\
    & = & b_{17}b_{14}b_{14}b_{13}b_{17}b_{15} \\
    & = & b_{13}b_{15}.
\end{eqnarray*}
This finishes the proof. $\Box$ \\ \\
Armed with the last two lemmas, we finish the proof of the reduction theorem via a series of technical lemmas. 
\begin{lemm} \mbox{}\\ 
(a) $[z,a,x,y]=[z,x,a,y]=[z,x,y,a]$. \\
(b) $[z,a,x,y,x]=[z,x,a,y,x]=[z,x,y,a,x]=[z,x,y,x,a]$.
\end{lemm}
{\bf Proof}\ \ (a) As  $[z,x,a]=[z,a,x]$, it suffices to show that $[z,x,a,y]=[z,x,y,a]$. Calculating in $\langle \bar{x},
\bar{a},\bar{b},\bar{y}\rangle=\langle x^{z},x,a,y\rangle$, we deduce that 
  $$[z,x,y,a]=[\bar{x}\bar{a},\bar{y},\bar{b}]=[[\bar{x},\bar{y}][\bar{x},\bar{y},\bar{a}][\bar{a},\bar{y}],\bar{b}]=e_{12},$$
and 
   $$[z,x,a,y]=[\bar{x},\bar{b},\bar{y}]=[e_{15},e_{18}]=e_{12}.$$
(b) By (a) the first three equalities hold,  and thus it only remains to show that $[z,x,a,y,x]=[z,x,y,x,a]$. Now $[z,x,a,y,x]=[e_{12},e_{19}]=e_{11}$. Also
\begin{eqnarray*}
      [z,x,y,x,a] & = & [\bar{x}\bar{a},\bar{y},\bar{a},\bar{b}] \\
\mbox{}    & = & [[\bar{x},\bar{y}][\bar{x},\bar{y},\bar{a}][\bar{a},\bar{y}],\bar{a},\bar{b}] \\
\mbox{}    & = & [\bar{x},\bar{y},\bar{a},\bar{b}] \\
\mbox{}   & = & [e_{14},e_{19},e_{20}] \\
\mbox{}   & = & [e_{13},e_{20}] \\
\mbox{} & =  & e_{11}.
\end{eqnarray*}
This finishes the proof. $\Box$
\begin{lemm} $[z,x,y,x,[z,a,y]]=[z,x,y,x,[z,a],y]=[z,a,x,y,x,[z,y]]=[z,a,x,y,x,z,y]$.
\end{lemm}
{\bf Proof}\ \ Notice that $[z,x,y,x]=[z,x,[y,x]]$ and $[z,a,x,y,x]=[z,a,x,[y,x]]$. Calculating in $\langle a,x,x^{y},z\rangle$, we
deduce that 
      $$[z,x,[y,x],[z,a]]=[z,a,x,[y,x],z]=b_{17}b_{14}.$$
Hence the 2nd and 4th terms are equal. The Hall-Witt identity implies that
\begin{eqnarray*}
  1 & = & [z,a,y,[z,x,y,x]]^{y}[y,[z,x,y,x],[z,a]]^{[z,x,y,x]}[z,x,y,x,[z,a],y]^{[z,a]} \\
    & = & [z,a,y,[z,x,y,x]][z,x,y,x,[z,a],y]^{[z,a]}.
\end{eqnarray*}
Thus
\begin{equation}
     [z,x,y,x,[z,a,y]]=[z,x,y,x,[z,a],y]^{[z,a]}.
\end{equation}
Notice that $[z,a,y]=[z,y,a]$. Using the Hall-Witt identity again, we get
\begin{eqnarray*}
   1 & = & [z,y,a,[z,x,y,x]]^{a}[a,[z,x,y,x],[z,y]]^{[z,x,y,x]}[z,x,y,x,[z,y],a]^{[z,y]} \\
     & = & [z,y,a,[z,x,y,x]]^{a}[a,[z,x,y,x],[z,y]],
\end{eqnarray*}
where we use the fact that $[z,x,y,x,a]$ commutes with $[z,x,y,x]$, which can for example be read from the presentation for 
$\langle x^{z},x,a,y\rangle$. Thus
\begin{equation}
    [z,x,y,x,[z,y,a]]=[z,x,y,x,a,[z,y]]^{a}.
\end{equation}
As $[z,x,[y,x],[z,a]]=b_{17}b_{14}$ commutes with $a$, and the same is true for $y$ and $[z,a]$, we see from (1) that $[z,x,y,x,[z,a,y]]$ commutes with
$a$. As this element also commutes with $z$, the conjugation by $a$ and $[z,a]$ in (1) and (2) can be dropped and the claimed equalities hold. $\Box$
\begin{lemm} $[z,x,y,x,[z,a,y,x]]=[z,x,y,x,[z,a,y],x]=
[z,x,y,x,[z,a],y,x]=[z,a,x,y,x,[z,y,x]].$
\end{lemm}  
{\bf Proof}\ \ It follows from Lemma 3.14 that the 2nd and 3rd terms are equal. We use the Hall-Witt identity to deduce that
 \begin{eqnarray*}
    1 & = & [z,a,y,x,[z,x,y,x]]^{x}[x,[z,x,y,x],[z,a,y,x]]^{[z,x,y,x]}[z,x,y,x,[z,a,y],x]^{[z,a,y]} \\
\mbox{} & = & [z,a,y,x,[z,x,y,x]][z,x,y,x,[z,a,y],x]^{[z,a,y]}.
\end{eqnarray*}
Thus
\begin{equation}
     [z,x,y,x,[z,a,y,x]]=[z,x,y,x,[z,a,y],x]^{[z,a,y]}.
\end{equation}
We also deduce that 
\begin{eqnarray*}
  1 & = & [z,y,x,a,[z,x,y,x]]^{a}[a,[z,x,y,x],[z,y,x]]^{[z,x,y,x]}[z,x,y,x,[z,y,x],a]^{[z,y,x]} \\
\mbox{} & = & [z,y,x,a,[z,x,y,x]]^{a}[a,[z,x,y,x],[z,y,x]],
\end{eqnarray*}
since $[a,[z,x,[y,x]]]$ commutes with $[z,x,[y,x]]$ in 
$\langle a,x,x^{y},z\rangle$. Thus
\begin{equation}
    [z,x,y,x,[z,y,x,a]]=[z,x,y,x,a,[z,y,x]]^{a}.
\end{equation}
From the proof of Lemma 3.14, we see that $[z,x,y,x,[z,a,y]]$ commutes with $a$. The same is true for $x$ and $[z,a,y]$ and thus it follows
from (3) that $[z,x,y,x,[z,a,y,x]]$ commutes with $a$. Calculating in the $\gamma$-group $\langle \bar{x},\bar{a},\bar{b},\bar{y}\rangle=
\langle y^{z},y,a,x\rangle$, we see that 
    $$[z,y,a,x,[z,y]]=[\bar{x},\bar{b},\bar{y},\bar{x}\bar{a}]=e_{11}e_{9}=1.$$
Thus $[z,x,y,x,[z,a,y,x]]$ commutes with $[z,y]$ and hence also with $[z,y,a]=[z,a,y]$. Therefore the conjugation by $a$ and $[z,a,y]$ can be dropped
and the lemma follows from (3) and (4) and Lemma 3.14. $\Box$
\begin{lemm}  $[z,x,y,x,[z,a,y,x,y]]=[z,x,y,x,[z,a,y,x],y]=[z,a,x,y,x,[z,y,x,y]]$.
\end{lemm}
{\bf Proof}\ \  The Hall-Witt identity implies that 
\begin{eqnarray*}
     1 & = & [z,a,y,x,y,[z,x,y,x]]^{y}[y,[z,x,y,x],[z,a,y,x]]^{[z,x,y,x]}[z,x,y,x,[z,a,y,x],y]^{[z,a,y,x]} \\
 \mbox{}      & = & [z,a,y,x,y,[z,x,y,x]][z,x,y,x,[z,a,y,x],y]^{[z,a,y,x]},
\end{eqnarray*}
so
  \begin{equation} 
[z,x,y,x,[z,a,y,x,y]]=[z,x,y,x,[z,a,y,x],y]^{[z,a,y,x]}.
\end{equation}
From the proof of Lemma 3.15, we see that $[z,x,y,x,[z,a,y,x]]$ commutes with $[z,a,y]$ and clearly $y$ commutes with $[z,a,y]$. It can be read from the presentation for $\langle y^{z},y,a,x\rangle$ that $[z,y,a,x]=[z,a,y,x]$ commutes with $[z,a,y]=[z,y,a]$. It follows that the RHS, and thus LHS, of (5) commutes
with $[z,a,y]$. Clearly the LHS also commutes with $x$. Hence it commutes with $[z,a,y,x]$ and 
\begin{equation}
    [z,x,y,x,[z,a,y,x,y]]=[z,x,y,x,[z,a,y,x],y].
\end{equation}
Also 
\begin{eqnarray*}
   1 & = & [z,y,x,y,a,[z,x,y,x]]^{a}[a,[z,x,y,x],[z,y,x,y]]^{[z,x,y,x]}[z,x,y,x,[z,y,x,y],a]^{[z,y,x,y]} \\
\mbox{} & = & [z,y,x,y,a,[z,x,y,x]]^{a}[a,[z,x,y,x],[z,y,x,y]].
\end{eqnarray*}
Thus
     $$[z,x,y,x,[z,y,x,y,a]]=[z,x,y,x,a,[z,y,x,y]]^{a}.$$
By the proof of Lemma 3.15, the RHS of (6) commutes with $a$. Hence
\begin{equation*}
    [z,x,y,x,[z,y,x,y,a]]=[z,x,y,x,a,[z,y,x,y]].
\end{equation*}
This finishes the proof. $\Box$
\begin{lemm} $[z,x,y,x,[z,a,y,x,y]]=[z,a,x,y,x,[z,y,x,y]]=1$. 
\end{lemm}
{\bf Proof}\ \ From Lemma 3.16, 
$$u:=[z,x,y,x,[z,a,y,x,y]]=[z,a,x,y,x,[z,y,x,y]].$$
Each of $[z,x,y,x]$ and $[z,a,y,x,y]$ commutes with each of $x,y$  and $[z,y]$.  That
$[z,a,y,x,y]=[z,y,a,x,y]$ commutes with $[z,y]$ can be read from the presentation for $\langle y^{z},y,a,x\rangle$. Thus the LHS commutes with $x,y,[z,y]$ and by symmetry the RHS commutes with $[z,x]$. Thus $u$ commutes with each of $x,y,[z,x]$ and $[z,y]$. We now show that
$u$ commutes with $a$: 
\begin{eqnarray*}
          u^{a} & = & [z,a,x,y,x,[z,y,x,y]^{a}] \\
\mbox{} & = & [z,a,x,y,x,[z,y,x,y][z,a,y,x,y]] \\
\mbox{}  & = & [z,a,x,y,x,[z,a,y,x,y][z,y,x,y]] \mbox{\ \ (working in }\langle y^{z},y,a,x\rangle)\\
\mbox{} & = & [z,a,x,y,x,[z,y,x,y]][z,a,x,y,x,[z,a,y,x,y]]^{[z,y,x,y]} \\
\mbox{} & = & [z,a,x,y,x,[z,y,x,y]] \\
\mbox{} & = & u.
\end{eqnarray*}
In the 3rd last identity we use the fact that $\langle a^{z},x,y\rangle$ is nilpotent of class at most $4$. \\ \\
We prove the following claim: if $v$ commutes with each of 
$a,x,y,[z,x]$ and $[z,y]$, then the same is true for $[v,[z,a]]$. 
To see this, notice first that as $[z,a]$ commutes with $[z,x]$,$[z,y]$ and $a$, the same is true for $[v,[z,a]]$. It follows that $v$ commutes with $[z,a,x]=[z,x,a]$ and thus
                     $$[v,[z,a]]^{x}=[v,[z,a][z,a,x]]=[v,[z,a,x][z,a]]=[v,[z,a]][v,[z,a,x]]^{[z,a]}=[v,[z,a]]$$
and $[v,[z,a]]$ commutes with $x$. Likewise $[v,[z,a]]$ commutes with $y$. Thus $[v,[z,a]]$ commutes with each of $a,x,y,[z,x]$ and $[z,y]$. \\ \\
Next, notice that $\langle [z,a]^{u},[z,a]\rangle\leq \langle a,z,z^{u}\rangle$ and thus $\langle [z,a]^{u},[z,a]\rangle$
is nilpotent, as every sandwich group of rank $3$ is nilpotent. 
Using this fact, let $m$ be the smallest positive integer such that
$[u,_{m}[z,a]]=1$. We claim that $m=1$. Notice that $[u,_{m-1}[z,a]]$ commutes 
with each of $x,y,[z,a],[z,x]$ and $[z,y]$. Proposition 3.3 implies that $[u,_{m-1}[z,a]]=1$. Hence $m=1$ and $u=1$. $\Box$
\begin{lemm} $v:=[z,x,y,x,[z,a,y,x]]=1$.
\end{lemm}
{\bf Proof}\ \ By Lemmas 3.17 and 3.16 we know that $v$ commutes with $y$. It also commutes with
$x$ (as $[a^{z},y,x,x]=1$). Observe that $[z,x,y,x,a]$ commutes with $[z,x,y,x]$ in $\langle x^{z},x,a,y\rangle$ and $[z,x,y,x,a]=[z,a,x,y,x]$ commutes
with $[z,a,y,x]$ in $\langle a^{z},x,y\rangle$. Hence
   $$v^{a}=[[z,x,y,x]^{a},[z,a,y,x]]=[[z,x,y,x,a][z,x,y,x],[z,a,y,x]]=[z,x,y,x,[z,a,y,x]]=v,$$
so $v$ commutes with $a$. It remains to show that $v$ commutes with $z$.  
Lemmas 3.14 and 3.15 imply that 
         $$v=[z,x,y,x,[z,a,y],x]=[z,a,x,y,x,[z,y],x].$$
We read from the presentation for $\langle x^{z},x,a,y\rangle$ that 
                $$[z,a,x,y,x]=[z,x,a,y,x]$$
 commutes with $[z,x]$ and by Lemma 3.10 this is also true for $[z,y]$.  Hence
      $$[z,a,x,y,x,[z,y]]=[z,x,y,x,[z,a,y]]$$
 commutes with $[z,x]$. Clearly $[z,x,y,x]$ and $[z,a,y]$ commute
with $z$. Thus
      $$v^{z}=[z,x,y,x,[z,a,y],x[x,z]]=[z,x,y,x,[z,a,y],[z,x]x]=v[z,x,y,x,[z,a,y],[z,x]]^{x}=v.$$
Thus $v\in Z(G)=1$. $\Box$
\begin{lemm} $v:=[z,x,y,x,[z,a,y]]=1$.
\end{lemm}
{\bf Proof}\ \ 
Lemmas 3.15 and 3.18 imply that $v$ commutes with $x$. It also commutes with $y$ and $z$. By the proof of Lemma 3.14,  it commutes with $a$. Thus $v\in Z(G)=1$. $\Box$
\begin{lemm} $v:=[z,x,y,x,[z,a]]=1$.
\end{lemm}
{\bf Proof}\ \ In $\langle a,x,x^{y},z\rangle$ we read that $v=b_{17}b_{14}$ commutes with $a,x$ and $z$. From Lemma 3.19 and 3.14, it commutes also with $y$. Thus $v\in Z(G)=1$. $\Box$ 
\begin{prop} $G$ is nilpotent.
\end{prop}
{\bf Proof}\ \ As $\langle x,y,z\rangle$ is nilpotent of class at most $4$, it follows from Lemma 3.20 and Proposition 3.3 that $[z,x,y,x]=1$. By symmetry in $x$ and $y$, we deduce that 
$[z,y,x,y]=1$. We saw earlier that $[z,x,[z,y]]=1$. It follows that $\langle x,y,z\rangle$ is nilpotent of class at most $3$. In particular, calculating in $\langle a,x,x^{y},z\rangle$, we see that 
            $$1=[z,x,[y,x],[z,a]]=b_{17}b_{14}.$$
Hence $b_{17}=b_{14}$. From the presentation for a $\beta$-group, 
we see that $b_{14}$ commutes with $b_{28}$. Hence $1=[b_{17},b_{28}]=b_{15}$. But $b_{15}=b_{13}$ so $b_{13}=1$. 
Also $b_{17}$ commutes with $a=b_{26}$ and $z=b_{18}$, and $x=b_{27}$. 
Recall that $b_{17}=[z,z^{a[x,y]}]$, so it is symmetric in $x$ and $y$; 
thus it commutes with $y$.  Hence $b_{17}\in Z(G)=\{1\}$. \\ \\
We now consider $[x,y,z]$. Calculating in $\langle a,x,x^{y},z\rangle$, we see that $[x,y,z,[z,a]]=b_{17}=1$. As $\langle x,y,z\rangle$ is nilpotent of class at most $3$,  Proposition 3.3 implies that $[x,y,z]=1$. But now
$[x,y]$ commutes with $a,x,y,z$ and thus $[x,y]=1$. Hence $G$ is a $\beta$-group and thus nilpotent. $\Box$ \\ \\
As a corollary to Theorem 3.2, we deduce the following.
\begin{theo} If $G=\langle a,x,y,z\rangle$ is a $\delta$-group, then $\langle z,z^{x},z^{y},a\rangle$ is nilpotent.
\end{theo}
{\bf Proof}\ \ Observe that  $\langle a_{1},x_{1},y_{1},z_{1}\rangle=\langle z,z^{x},z^{y},a\rangle $ is also a $\delta$-group. Now 
        $$[z_{1},a_{1},x_{1},z_{1}]=[a,z,z^{x},a]=[a^{zz^{x}},a]=[z,x,a,a]=[x^{z},a,a]=1.$$
Similarly we see that $[z_{1},a_{1},y_{1},z_{1}]=1$. It follows from Theorem 3.2 that $\langle a_{1},x_{1},y_{1},z_{1}\rangle$ is nilpotent. $\Box$. 
					
\subsubsection{Sandwich groups whose commutativity graph has $1$ edge}
Let $G$ be the largest sandwich group generated by 4 involutions
where precisely one pair commutes.
We used the {\sc Magma} implementation of the
nilpotent quotient algorithm to determine that
the largest nilpotent quotient of 
$G$ has order $2^{235}$ and class~11.

\subsection{Rank 5 sandwich groups generated by involutions}
Let $G$ be a residually nilpotent sandwich group of rank $5$ 
generated by involutions.
It remains an open question whether $G$ is finite.
As a step towards its resolution, we report the following.

\begin{theo}
A residually nilpotent sandwich group of rank $5$ generated by involutions
is finite if there are least three commuting pairs among its generators.
\end{theo}
We consider all 10 commutativity graphs for such a group $G$ that have 
$3$ or $4$ edges. 
In Table \ref{summary-bounds}
we list the graphs  and upper 
bounds for the class and orders of 
the largest 2-quotient of each group.

We proved these claims computationally by
studying each of the 10 groups using our implementation
in {\sc Magma} of the $p$-quotient algorithm \cite{pQ}.
Let $Q$ be the group generated by the involutions $a,b,c,d,e$
which satisfy one of these commutativity graphs.
It required too much CPU time to impose the sandwich condition
directly. Instead, we added up to 2500 random instances of
this condition as explicit relations to the presentation for $Q$,
and constructed its largest finite $2$-quotient. The most expensive of these computations took 7 days of CPU time using {\sc Magma} 2.27-3 on a computer with a 2.6 GHz processor. 
Computational evidence suggests that the largest 2-quotient
of $G$ and $Q$ coincide.

\begin{table}[ht]
$$\begin{array}{l|c|r}
\mbox{Commutativity graph} & \mbox{Class bound} & \log_{2}|G| \\  \hline
a-b,a-c,a-d,a-e & 13 & 777 \\ \hline 
a-b-c-d-e & 17 & 2643 \\ \hline
a-b-c-d-a & 16 & 2831 \\ \hline
a-b-c-a, d-e & 16 & 3145 \\ \hline
a-b-c-a-d & 16 & 3324 \\ \hline 
a-b, a-c, a-d-e & 16 & 2636 \\ \hline 
a-b-c-d  & 17 & 10354 \\ \hline 
a-b-d-a  & 16 & 12598 \\ \hline 
a-b, a-c, a-d & 16 & 10906 \\ \hline 
a-b-c, d-e & 16 & 9987
\end{array}$$
\caption{Class and order bounds for the 2-quotients}\label{summary-bounds}
\end{table}

\section{Global nilpotence question for $\langle a\rangle^{G}$ in locally finite $2$-groups}

In  this section we give an example from \cite{GGM} of a locally finite $2$-group $G$ with a left $3$-Engel element $a$ such that
$\langle a\rangle^{G}$ is not nilpotent. The construction is based on a Lie algebra given in \cite{Trau2}. This result was generalised in \cite{GA} to an infinite family of examples.  

\subsection{The Lie algebra}

Let $\mathbb{F}$ be the field of order 2 and consider the 4-dimensional vector space $V=\mathbb{F} x+\mathbb{F} u+\mathbb{F} v+ \mathbb{F} w $. We equip $L$ with a binary product where
$$u \cdot v = u, \quad v \cdot w = w, \quad w \cdot u = v, \quad u \cdot x = 0, \quad v \cdot x = 0, \quad w \cdot x = u.$$
We then extend the product linearly on $V$. Observe that $V$ is a Lie algebra with a trivial center and $W={\mathbb F}u+{\mathbb F}v+{\mathbb F}w$ is a simple ideal. \\ \\
Let $E=\langle \mbox{ad}(x),\mbox{ad}(u),\mbox{ad}(v),\mbox{ad}(w)\rangle\leq \mbox{End}(V)$ be the associative enveloping algebra of $V$. Now $E$ is $12$-dimensional with basis
$$\begin{array}{llll}
     e_{1}=\mbox{ad}(w), & e_{2}=\mbox{ad}(w)^{2}, & e_{3}=\mbox{ad}(w)^{3},  & e_{4}=\mbox{ad}(v), \\
   e_{5}=\mbox{ad}(v)\mbox{ad}(w), & e_{6}=\mbox{ad}(v)\mbox{ad}(w)^{2},  &  e_{7}=\mbox{ad}(u), & e_{8}=\mbox{ad}(u)\mbox{ad}(w), \\
   e_{9}=\mbox{ad}(u)\mbox{ad}(w)^{2}, & e_{10}=\mbox{ad}(x)\mbox{ad}(v), & e_{11}=\mbox{ad}(x)\mbox{ad}(w), & e_{12}=\mbox{ad}(x)\mbox{ad}(w)^{2}.
\end{array}$$
We construct a certain locally nilpotent Lie algebra over ${\mathbb F}$ of countably infinite dimension. For ease of notation, we introduce the following modified union of subsets of ${\mathbb N}$:
\[A  \sqcup B=\begin{cases} A \cup B, & \text{ if } A \cap B = \emptyset \\
\emptyset & \text{ otherwise}.
\end{cases}\]
For each non-empty subset $A$ of ${\mathbb N}$, let $W_{A}$ be a copy of $W$. That is, $W_{A}=\{z_{A}:z\in W\}$ with
addition $z_{A}+t_{A}=(z+t)_{A}$. We  take the direct sum 
     $$W^{*}=\bigoplus_{\emptyset\not =A\subseteq {\mathbb N}}W_{A}.$$
We view $W^{*}$ as a Lie algebra by defining a multiplication 
         $$z_{A}\cdot t_{B}=(zt)_{A\sqcup B},$$
for $z_{A}\in W_{A}$ and $t_{B}\in W_{B}$, and then extend this product linearly on $W^{*}$. The interpretation here is that $z_{\emptyset}=0$. Finally, we extend this to a semidirect product with ${\mathbb F}x$ 
             $$V^{*}=W^{*}\oplus {\mathbb F}x$$
induced from the action $z_{A}\cdot x=(zx)_{A}$.  \\ \\
Notice that $V^{*}$ has basis 
       $$\{x\}\cup\{u_{A},v_{A},w_{A}:\emptyset\not =A\subseteq {\mathbb N}\}$$
and 
        $$u_{A}\cdot u_{B}=v_{A}\cdot v_{B}=w_{A}\cdot w_{B}=0,$$
     $$u_{A}\cdot x=0,\ v_{A}\cdot x=0,\ w_{A}\cdot x=u_{A}$$
and
    $$u_{A}\cdot v_{B}=u_{A\sqcup B},\ v_{A}\cdot w_{B}=w_{A\sqcup B},\ w_{A}\cdot u_{B}=v_{A\sqcup B}.$$
Every finitely generated subalgebra of $V^{*}$ is contained in some 
$$S=\langle x, u_{A_{1}}, \ldots ,
u_{A_{r}}, v_{B_{1}},\ldots ,v_{B_{s}},w_{C_{1}},\ldots ,w_{C_{t}}\rangle.$$
 Since $zxx=0$ for all $z\in V^{*}$, it follows that $S$ is nilpotent
of class at most $2(r+s+t)$. Hence $V^{*}$ is locally nilpotent. \\ \\
We now construct a group $G\leq \mbox{GL}(V^{*})$ containing $a:=1+\mbox{ad}(x)$ where
$a$ is a left $3$-Engel element of $G$ but $\langle a\rangle^{G}$ is not nilpotent. 
Let $y$ be one of $x,u_{A},v_{A},w_{A}$. Observe that $\mbox{ad}(y)^{2}=0$ so 
      $$(1+\mbox{ad}(y))^{2}=1+2\mbox{ad}(y)+\mbox{ad}(y)^{2}=1.$$
Thus $1+\mbox{ad}(y)$ is an involution in $\mbox{GL}(V^{*})$. The subgroups 
        $${\mathcal U}=\langle 1+\mbox{ad}(u_{A}):\,A\subseteq {\mathbb N}\rangle,\,
				{\mathcal V}=\langle 1+\mbox{ad}(v_{A}):\,A\subseteq {\mathbb N}\rangle,\ 
				{\mathcal W}=\langle 1+\mbox{ad}(w_{A}):\,A\subseteq {\mathbb N}\rangle
				$$
are elementary abelian of countably infinite rank. 
%
Analysis 
of $G:=\langle a, {\mathcal U},{\mathcal V},{\mathcal W}\rangle$
establishes the following. 
\begin{theo} The element $a$ is a left $3$-Engel element of $G$. However $\langle a\rangle^{G}$ is not nilpotent. 
\end{theo}

\section{Global nilpotence question for $\langle x\rangle^{G}$ in locally finite $p$-groups, for odd $p$}

In Section 4, following \cite{GA, GGM}, 
we gave an example of a locally finite $2$-group $G$ with a left $3$-Engel element $a$ such that $\langle a\rangle^{G}$ is not nilpotent. 
We now 
provide such an example for locally finite $p$-groups where $p$ is any odd prime \cite{GAM}. The odd case is more involved and the construction quite different from the $p=2$ case. We first describe the pair $(G,x)$ that will provide our example and we show directly that $x$ is a left $3$-Engel element of $G$. The description of $G$ is not as transparent as in Section 4.  To show that $\langle x\rangle^{G}$ is not nilpotent requires more work. We construct first a pair $(L,z)$ where $L$ is a Lie algebra over ${\mathbb F}_{p}$, the field of $p$ elements, and $\mbox{Id}(z)$ is not nilpotent. The pair $(L,z)$ can be seen as the Lie algebra analogue of our group construction. We then build a group $H$ within $\mbox{End}(L)$ containing $1+\mbox{ad}(z)$ where $(H, 1+\mbox{ad}(z))$ is a homomorphic image of $(G,x)$. Since $(1+\mbox{ad}(z))^{H}$ is not nilpotent,
$\langle x\rangle^{G}$ is not nilpotent.

\subsection{The group $G$ and Lie algebra $L$}
Let $x, a_{1}, a_{2},\ldots $ be an infinite list of group variables. Recall that a simple commutator in $x, a_{1}, a_{2}, \ldots $ is a group word defined recursively as follows: $x, a_{1}, a_{2}, \ldots $ are simple commutators; if $u$ and $v$ are simple commutators, then $[u,v]$ is a simple commutator.  A simple commutator $s$ has multi-weight $(m,e_{1},e_{2},\ldots )$ in $x, a_{1}, a_{2}, \ldots $, if $x$ occurs $m$ times and $a_{i}$ occurs $e_{i}$ times in $s$. The weight of $s$ is $m+e_{1}+e_{2}+\cdots $. The following definition is critical.\\ \\
{\bf Definition}. Let $s$ be a simple commutator of multi-weight $(m,e_{1},e_{2},\ldots )$ in  $x,a_{1},a_{2},\ldots$. The {\it type} of $s$ is $t(s)=e_{1}+e_{2}+\ldots -2m$.\\ \\
{\bf Remark}. If $u,v$ are simple commutators in $x,a_{1},a_{2},\ldots$, then $t([u,v])=t(u)+t(v)$. In particular, $t([u,a_{j}])=t(u)+1$
and $t([u,x])=t(u)-2$. \\ \\
For a fixed odd prime $p$, let $G= \langle x, a_1, a_{2}, \dots \rangle$ be the largest group satisfying the following conditions:
\begin{enumerate}
\item[(1)] $\langle a_i \rangle^G$ is abelian for all $i\geq 1$; 
\item[(2)] $\langle x \rangle^G$ is metabelian; 
\item[(3)] $x^{p}=a_{1}^{p}=a_{2}^{p}=\cdots =1$; 
\item[(4)] if $s \neq x$ is a simple commutator in 
$x,a_{1},a_{2},\ldots$ and 
$|t(s)| \geq 2$, 
then $s=1$.
\end{enumerate}
%
It is not difficult to see that $G$ is a locally finite $p$-group. From (4), it is clear that $s=1$ unless it is of the form 
 $[x,a_{j_{1}},a_{j_{2}},a_{j_{3}},x,a_{j_{4}},a_{j_{5}},  \dots, x,a_{j_{2m}},a_{j_{2m+1}}]$. \\ \\
In \cite{GAM} it was proved directly that $x$ is a left $3$-Engel element of $G$. To show that $\langle x\rangle^{G}$ is not nilpotent, it suffices to show that the special commutators above, which are not trivial by (4), are non-trivial. The structure of $G$ is not transparent enough for us to prove this directly. Instead we look at the analogous Lie algebra setting. \\ \\
We first consider  the  largest Lie algebra $F=\langle z, c_1, c_2, \dots \rangle$ over $\mathbb{F}_p$ such that: 
\begin{enumerate}
    \item $\mbox{Id}(c_i)$ is abelian for $i=1,2, \dots$;
    \item $\mbox{Id}(z)$ is metabelian.
\end{enumerate}
%
Let $B$ be the following basis for $\mbox{Id}_F(z)$: 
\begin{equation*}
[[z, c_{I_1}],[z, c_{I_2}], \dots, [z, c_{I_m}]]
\end{equation*}
where $m \geq 1$, $I_1, \dots, I_m$ are pairwise disjoint and $[z, c_{I_1}] > [z, c_{I_2}]\leq \dots \leq [z, c_{I_m}]$. \\ \\
We define the type of a Lie commutator by analogy with that of a group commutator. \\ \\
{\bf Definition}. Let $s$ be a simple commutator of multi-weight $(m,e_1, \dots, e_r)$ in $z, c_1, \dots, c_r$. The {\it type} of $s$ is $t(s)=e_1+\dots +e_r-2m$.\\ \\
{\bf Remark}. If $c$ and $d$ are simple commutators, then $t([c,d])=t(c)+t(d)$. In particular,
$t([c,x])=t(c)-2$ and $t([c,a_{1}])=t(c)+1$.\\ \\
We construct a Lie algebra $L$ that is a quotient of $F$ by a certain multi-homogeneous ideal $J$. Since $F$ is multi-graded, $L$ is also multi-graded.  
Now $J$ is the smallest ideal containing the following: 
\begin{enumerate}
    \item[3.] $c \in B$ is in $J$ if $c \neq z$ is a commutator of type having absolute value greater than 1;
    \item[4.] $c \in B$ is in $J$ if one of $I_1,I_3, \dots, I_m$ has size greater than 2.
\end{enumerate}
By careful analysis, we determine the structure of $J$ and prove that $\mbox{Id}_L(z)$ is not nilpotent. \\ \\
From $L$, we obtain a group 
$H=\langle 1+\mbox{ad}(z),1+\mbox{ad}(c_1),1+\mbox{ad}(c_2), \dots \rangle$. 
We deduce that $H$ is a homomorphic image of $G$ by showing that 
relations (1)--(4) of $G$ hold in $H$ where
$x, a_{1}, a_{2},\dots $ are replaced by 
$1+\mbox{ad}(z), 1+\mbox{ad}(c_1), 1+\mbox{ad}(c_2),\ldots $. 
Finally, we show that $\langle 1+\mbox{ad}(z)\rangle^{H}$ is non-nilpotent, 
so $\langle x\rangle^{G}$ is non-nilpotent. 

\begin{theo} The element $x$ is a left $3$-Engel element of $G$. However $\langle x\rangle^{G}$ is not nilpotent. 
\end{theo}
%


\medskip
\noindent
Hadjievangelou and Traustason: Department of Mathematical Sciences, 
University of Bath, Claverton Down, Bath BA2 7AY, UK

\medskip
\noindent
Longobardi, Maj and Monetta: Universit\`a Degli Studi di Salerno, Via Giovanni Paolo II, 132 - 84084 Fisciano (SA), Italy

\medskip
\noindent
O'Brien, Department of Mathematics, University of Auckland, Private
Bag 92019, Auckland, New Zealand


\vspace{0.3cm}
\begin{thebibliography}{99}
%
\bibitem{Ab1} A. Abdollahi, Left $3$-Engel elements in groups, 
{\it J. Pure Appl. Algebra}
{\bf 188} (2004), 1-6.
%
\bibitem{Ab2} A. Abdollahi and H. Khosravi, On the right and left 4-Engel elements, {\it Comm. Algebra} {\bf 38} (2010), 993-943.
%
\bibitem{Baer} R. Baer, Engelsche Elemente Noetherscher Gruppen, {\it Math. Ann.} {\bf 133} (1957), 256-270. 
%
\bibitem{MAGMA} W. Bosma, J. Cannon and C. Playoust,
The Magma Algebra System I: The User Language, {\it J. Symbolic Comput.}, {\bf
24} (1997), 235-265.
%
\bibitem{Burn} W. Burnside, On an unsettled question in the theory of 
discontinuous groups, {\it Quart.\ J. Pure Appl.\ Math.} {\bf 37} (1901), 230-238.
%
%
%
\bibitem{Gru} K. W. Gruenberg, The Engel elements of a soluble group, {\it Illinois J. Math.} {\bf 3} (1959), 151-169.
%
\bibitem{GA1} A. Hadjievangelou and G. Traustason, Sandwich groups and (strong) left $3$-Engel elements in groups, submitted.
%
\bibitem{GA} A. Hadjievangelou and G. Traustason, Left 3-Engel elements in locally finite 2-groups, {\it Comm.  Algebra} {\bf 49} (2021), 4869-4882.
%
\bibitem{GAM} A. Hadjievangelou, M. Noce and G. Traustason, Locally finite $p$-groups with a left 3-Engel element whose normal closure is not nilpotent, {\it Internat. J. Algebra Comput.} {\bf 31} (2021), 135-160. 
%
\bibitem{GAP} The GAP Group, GAP -- Groups, Algorithms, and Programming, Version 4.12.1; 2022.
%
\bibitem{Hein1} H. Heineken, Eine Bemerkung \"{u}ber engelshe Elemente, {\it Arch. Math.} {\bf 11} (1960), 321.
%
\bibitem{CGT-Handbook}
D. F.\ Holt, B. Eick, and E. A.\ O'Brien,
\newblock {\em Handbook of computational group theory},
\newblock Chapman and Hall/CRC, London, 2005
%
\bibitem{Ivan} S. V. Ivanov, The free Burnside groups of sufficiently large
exponents, {\it Internat. J. Algebra Comput.} {\bf 4} (1994), 1-308.
%
\bibitem{Ko1} A. I. Kostrikin, {\it Around Burnside}, Moscow, Nauka (1986).
%
\bibitem{Ko2} A. I. Kostrikin, The Burnside Problem, {\it Izv. Akad Nauk SSSR, Ser. Mat.}, {\bf 23} (1959), 3-34. 
%
\bibitem{Kos} A. I. Kostrikin and E. I. Zel'manov, A Theorem on Sandwich Algebras, {\it Trudy Mat. Inst. Steklov}, {\bf 183} (1988), 142-149.
%
\bibitem{Jab} E. Jabara and G. Traustason, Left $3$-Engel elements of odd order in groups, {\it Proc. Amer. Math. Soc.} {\bf 147} (2019), 1921-1927.
%
\bibitem{Lys} I. G. Lysenok, Infinite Burnside groups of even period,
{\it Izv. Ross. Akad. Nauk. Ser. Math.} {\bf 60} (1996), 453-654.
%
\bibitem{New} M. Newell, On right-Engel elements of length three, {\it Proc. Roy. Irish Acad.} {\bf 96A}, No.1 (1996), 17-24. 
%
\bibitem{pQ} M. F.\ Newman and E. A. O'Brien, Application of computers to questions like those of Burnside, {II},
{\it Internat.\ J.\ Algebra Comput.} {\bf 6}, 593-605, 1996.
%
\bibitem{NQ} W. Nickel,
Computing nilpotent quotients of finitely presented groups.
Geometric and computational perspectives on infinite groups
(Minneapolis, MN and New Brunswick, NJ, 1994), 175-191,
DIMACS Ser. Discrete Math. Theoret. Comput. Sci., 25,
Amer. Math. Soc., Providence, RI, 1996.
%
\bibitem{GGM} M. Noce, G. Tracey and G. Traustason, A left $3$-Engel element whose normal closure is not nilpotent, {\it J.  Pure Appl. Algebra} {\bf 224}, (2020), 1092-1101.
%
\bibitem{San} I. N. Sanov, Solution of Burnside's Problem for Exponent Four, Leningrad State Univ., {\it Ann. Maths. Ser. } {\bf 10} (1940), 166-170.  
%
\bibitem{Trac} G. Tracey and G. Traustason, Left $3$-Engel elements in groups of exponent $60$, {\it Internat. J. Algebra  Comput.} {\bf 28} (2018), 673-695. 
%
\bibitem{Trau2} G. Traustason, Engel Lie-algebras, {\it Quart. J. Math. Oxford.} Ser. (2), {\bf 44} (1993), 355-384. 
%
%
\bibitem{Trau} G. Traustason, Left $3$-Engel elements in groups of exponent $5$, {\it J. Algebra}, {\bf 414} (2014), 41-71.
%
\bibitem{Vaug} M. Vaughan-Lee,  {\it The Restricted Burnside Problem}, Clarendon Press Oxford (2nd edition) (1993).
%
\bibitem{Vaug1} M. Vaughan-Lee, Engel-$4$ groups of exponent $5$, {\it Proc. Lond. Math. Soc.}, {\bf 74} (1997), 306-334.
%
\bibitem{zc} E. I. Zel'manov, The solution of the restricted Burnside problem for groups of odd exponent, {\it Math. USSR Izvestia} {\bf 36} (1991), 41-60.
%
\bibitem{zd} E. I. Zel'manov, The solution of the restricted Burnside problem for $2$-groups, {\it Mat. Sb.} {\bf 182} (1991), 568-592.

\end{thebibliography}
\end{document}